\documentclass[]{article}

\usepackage[margin=30truemm]{geometry}

\usepackage[whole]{bxcjkjatype}

\usepackage[cmex10]{amsmath}
\usepackage{amssymb}
\usepackage{amsthm}

\usepackage{cases} 

\usepackage{fixltx2e}

\usepackage[subrefformat=parens]{subcaption}
\usepackage{graphicx}
\usepackage{epstopdf}

\usepackage{xcolor}

\usepackage{cite}

\usepackage{algorithm}
\usepackage{algorithmic}

\newtheorem{definition}{Definition}
\newtheorem{theorem}{Theorem}
\newtheorem{lemma}{Lemma}
\newtheorem{proposition}{Proposition}
\newtheorem{corollary}{Corollary}

\newtheorem{assumption}{Assumption}
\newtheorem{remark}{Remark} 
\newtheorem{example}{Example}

\newcommand{\argmin}{\operatornamewithlimits{\mathrm{arg~min}}}

\definecolor{myGreen}{rgb}{ 0, 0.7, 0.3 }
\definecolor{myBlue}{rgb}{ 0, 0.4, 1 }
\definecolor{myPurple}{rgb}{ 0.7, 0, 0.3 }
\definecolor{myGray}{gray}{ 0.7 }
\definecolor{myFigGray}{gray}{ 0.5 }

\usepackage{pgfplots}
\pgfplotsset{compat=newest}
\usepackage{tikz}
\usetikzlibrary{positioning}
\usetikzlibrary{arrows}

\newcommand{\Mymemo}[1]{}

\usepackage[overlay]{textpos}

\usepackage{pgfplots}
\pgfplotsset{compat=newest}
\usepackage{tikz}
\usetikzlibrary{positioning}
\usetikzlibrary{arrows}

\usepackage{natbib}

\begin{document}
	
	\title{Weighted Stochastic Riccati Equations for Generalization of Linear Optimal Control%
		\thanks{%
			This paper is submitted to a journal for possible publication. The copyright of this paper may be transferred without notice, after which this version may no longer be accessible.
			This work was partly supported by JSPS KAKENHI Grant Number JP18K04222.
			The material herein was presented in part at the 2016 American Control Conference \citep{ItoACC16}.
			We would like to thank Editage (www.editage.jp) for the English language editing.
	}} 
	
	\author{Yuji Ito\thanks{Yuji Ito is the corresponding author and with Toyota Central R\&D Labs., Inc., 41-1 Yokomichi, Nagakute-shi, Aichi 480-1192, Japan	(e-mail: ito-yuji@mosk.tytlabs.co.jp).}
		\and
		Kenji Fujimoto\thanks{Kenji Fujimoto is with the Department of Aeronautics and Astronautics, Graduate School of Engineering, Kyoto University, Kyotodaigakukatsura, Nishikyo-ku, Kyoto-shi, Kyoto 615-8540, Japan (e-mail: k.fujimoto@ieee.org).}
		\and
		Yukihiro Tadokoro\thanks{Yukihiro Tadokoro is with Toyota Central R\&D Labs., Inc., 41-1 Yokomichi, Nagakute-shi, Aichi 480-1192, Japan	(e-mail: y.tadokoro@ieee.org).}		
	}  
	
	\date{}
	
	\maketitle

\begin{abstract}
This paper presents weighted stochastic Riccati (WSR) equations for designing multiple types of optimal controllers for linear stochastic systems.
The stochastic system matrices are independent and identically distributed (i.i.d.) to represent uncertainty and noise in the systems.
However, it is difficult to design multiple types of controllers for systems with i.i.d. matrices while the stochasticity can invoke unpredictable control results.
A critical limitation of such i.i.d. systems is that Riccati-like algebraic equations cannot be applied to complex controller design.
To overcome this limitation, the proposed WSR equations employ a weighted expectation of stochastic algebraic equations.
The weighted expectation is calculated using a weight function designed to handle statistical properties of the control policy. 
Solutions to the WSR equations provide multiple policies depending on the weight function, which contain the deterministic optimal, stochastic optimal, and risk-sensitive linear (RSL) control.
This study presents two approaches to solve the WSR equations efficiently: calculating WSR difference equations iteratively and employing Newton's method.
Moreover, designing the weight function yields a novel controller termed the robust RSL controller that has both a risk-sensitive policy and robustness to randomness occurring in stochastic control design.
\end{abstract}

\newcommand*{\SymColor}[1]{\textcolor{red}{#1}}
\renewcommand*{\SymColor}[1]{#1}

\newcommand*{\KoningTransition}{\SymColor{(i)}}
\newcommand*{\KoningStable}{\SymColor{(ii)}}
\newcommand*{\KoningStabilizable}{\SymColor{(iii)}}
\newcommand*{\KoningTransStable}{\SymColor{(iv)}}
\newcommand*{\KoningPiConverge}{\SymColor{(v)}}
\newcommand*{\KoningIHopt}{\SymColor{(vi)}}
\newcommand*{\KoningUnique}{\SymColor{(vii)}}

\newcommand*{\StatementImplicitFunction}{\SymColor{(S1)}}
\newcommand*{\StatementBallforImplicitFunction}{\SymColor{(S2)}}
\newcommand*{\StatementSetStandardAss}{\SymColor{(S3)}}
\newcommand*{\StatementInitialGuess}{\SymColor{(S4)}}

\newcommand*{\StandardAssA}{\SymColor{(i)}}
\newcommand*{\StandardAssB}{\SymColor{(ii)}}
\newcommand*{\StandardAssC}{\SymColor{(iii)}}

\newcommand*{\AssWeightMSstabilizable}{\SymColor{(i)}}
\newcommand*{\AssWeightZeroSolContained}{\SymColor{(ii)}}
\newcommand*{\AssWeightSmooth}{\SymColor{(iii)}}

\newcommand*{\El}[2]{\SymColor{[}{#1}\SymColor{]}_{#2}}
\newcommand*{\VEC}[2][]{\SymColor{\mathrm{vec}#1(}   #2 \SymColor{#1)}}
\newcommand*{\VECH}[2][]{\SymColor{\mathrm{vech}#1(} #2 \SymColor{#1)}}

\newcommand*{\FrobeniusNorm}[1]{\SymColor{\|}{#1}\SymColor{\|_{\mathrm{F}}}} 

\newcommand*{\Distance}[2]{ \SymColor{d(} #1 \SymColor{,} #2 \SymColor{)}  }

\newcommand*{\DIAG}[1]{\SymColor{\mathrm{diag}(} #1\SymColor{)}}

\newcommand*{\MyTRANSPO}{\SymColor{\top}}

\newcommand*{\Elimi}[1]{\SymColor{\mathcal{C}} \SymColor{(}    #1\SymColor{)}}
\newcommand*{\DupMat}[1]{\SymColor{\mathcal{D}}_{\SymColor{#1}}}
\newcommand*{\EliMat}[1]{\SymColor{\mathcal{L}}_{\SymColor{#1}}}

\newcommand*{\wildcard}{  \SymColor{\bullet}  }

\newcommand*{\Expect}[1]{\SymColor{\mathrm{E}}_{#1}}
\newcommand*{\Covariance}[1]{\SymColor{\mathrm{Cov}}_{#1}}

\newcommand*{\WE}[4]{\SymColor{\mathcal{E}}_{#1}^{(#2,#3,#4)}  }
\newcommand*{\biasE}[1]{\SymColor{\mathrm{E}}_{#1}^{\SymColor{\mathrm{q}}}  }

\newcommand*{\PDF}[1]{\SymColor{f_{\mathrm{p}}(}#1\SymColor{)}}
\newcommand*{\biasPDF}[1]{\SymColor{f_{\mathrm{q}}(}#1\SymColor{)}}


\newcommand*{\Identity}[1]{\SymColor{\boldsymbol{I}}_{ #1 }}

\newcommand*{\SetSymMat}[1]{\SymColor{\mathbb{R}}_{\mathrm{sym}}^{#1}}

\newcommand*{\NotationVec}{\SymColor{\boldsymbol{v}}}
\newcommand*{\NotationMat}{\SymColor{\boldsymbol{C}}}
\newcommand*{\NotationSquMat}{\SymColor{\boldsymbol{D}}}
\newcommand*{\NotationSymMat}{\SymColor{\boldsymbol{S}}}

\newcommand*{\NotationSymAMat}{\NotationSymMat_{\NotationA}}
\newcommand*{\NotationSymBMat}{\NotationSymMat_{\NotationB}}
\newcommand*{\NotationSto}{\SymColor{\boldsymbol{\Lambda}}}
\newcommand*{\NotationFunc}{\SymColor{\boldsymbol{\phi}}}
\newcommand*{\IDNotation}{\SymColor{i}}
\newcommand*{\IDbNotation}{\SymColor{j}}
\newcommand*{\DimANotation}{\SymColor{a}}
\newcommand*{\DimBNotation}{\SymColor{b}}
\newcommand*{\NotationA}{\SymColor{1}}
\newcommand*{\NotationB}{\SymColor{2}}

\newcommand*{\IDiteA}{\SymColor{\ell}}

\newcommand*{\IDEl}{\SymColor{i}}
\newcommand*{\IDbEl}{\SymColor{j}}

\newcommand*{\DimX}{\SymColor{n}}
\newcommand*{\DimU}{\SymColor{m}}
\newcommand*{\DimexX}{\SymColor{\tilde{n}}}

\newcommand*{\DomALLFBgain}{\SymColor{\overline{\mathbb{S}}_{\FBgain}}}
\newcommand*{\DomALLVmat}{\SymColor{\overline{\mathbb{S}}_{\Vmat}}}
\newcommand*{\DomALLSolPair}{\SymColor{\overline{\mathbb{S}}_{\SolPair}}}

\newcommand*{\DomStoParam}{\SymColor{\mathbb{S}_{\StoParam}}}
\newcommand*{\DomFBgain}{\SymColor{\mathbb{S}_{\FBgain}}}
\newcommand*{\DomVmat}{\SymColor{\mathbb{S}_{\Vmat}}}
\newcommand*{\DomSolPair}{\SymColor{\mathbb{S}_{\SolPair}}}

\newcommand*{\DomSaSolPair}{\SymColor{\mathbb{S}_{\SolPair}^{\prime}}}
\newcommand*{\BallSolPair}[2]{\SymColor{{\mathbb{B}}_{\SolPair}^{\ast}} \SymColor{(} #1 \SymColor{,} #2 \SymColor{)}} 

\newcommand*{\DomBSolPair}{\SymColor{\mathbb{S}_{\SolPair,2}}}
\newcommand*{\DomNeighborSolPair}{\SymColor{\mathbb{B}_{\SolPair}}}
\newcommand*{\DomNeighborriskParam}{\SymColor{\mathbb{B}_{\riskParam}}}

\newcommand*{\BallNonsingular}{   \SymColor{\mathbb{B}_{\ast}}    }


\newcommand*{\MyT}{\SymColor{t}}
\newcommand*{\MybT}{\SymColor{s}}
\newcommand*{\State}{\SymColor{\boldsymbol{x}}}
\newcommand*{\Input}{\SymColor{\boldsymbol{u}}}
\newcommand*{\DriftMat}{\SymColor{\boldsymbol{A}}}
\newcommand*{\InMat}{\SymColor{\boldsymbol{B}}}
\newcommand*{\StoParam}{\SymColor{\boldsymbol{\Lambda}}}
\newcommand*{\FBgain}{\SymColor{\boldsymbol{L}}}
\newcommand*{\DescriptionFBMat}[1]{{(\DriftMat - \InMat #1)}}

\newcommand*{\RefX}{\SymColor{\tilde{\State}}}

\newcommand*{\riskParam}{\SymColor{\theta}}

\newcommand*{\terminalT}{\SymColor{T}}

\newcommand*{\xCostMat}{\SymColor{\boldsymbol{Q}}}
\newcommand*{\uCostMat}{\SymColor{\boldsymbol{R}}}

\newcommand*{\FHWcostJ}[4]{\SymColor{J}_{#4}(#1,#2;#3)}
\newcommand*{\IHWcostJ}[3]{\SymColor{J}_{\infty}(#1,#2;#3)}
\newcommand*{\funccostJ}{     \SymColor{J}   }
\newcommand*{\estfunccostJ}{        \SymColor{\widetilde{J}}   }  
\newcommand*{\mincostJ}{     \SymColor{J_{\ast}}    }
\newcommand*{\predcostJ}{\SymColor{J_{\mathrm{E}}}}

\newcommand*{\funcInput}{       \SymColor{{\boldsymbol{u}}}            }     
\newcommand*{\estfuncInput}{       \SymColor{\widetilde{\boldsymbol{u}} }           } 
\newcommand*{\optFBcontroller}{     \SymColor{\boldsymbol{u}_{\ast}}           }

\newcommand*{\dWeight}{\SymColor{w_{\mathrm{d}}}}
\newcommand*{\Weight}[1]{\SymColor{w ( \StoParam_{#1} ; \riskParam )}}

\newcommand*{\paramWeightA}{\SymColor{\alpha}}
\newcommand*{\paramWeightB}{\SymColor{\beta}}

\newcommand*{\WmsWeight}[1]{   \SymColor{W}_{#1}\SymColor{(} \StoParam_{0:#1}  \SymColor{)}       }


\newcommand*{\WSRDVmat}[3]{       \SymColor{\boldsymbol{F}}         \SymColor{(} #1, #2 ; #3 \SymColor{)}       }
\newcommand*{\WSRDFBgain}[3]{     \SymColor{\boldsymbol{G}}         \SymColor{(} #1, #2 ; #3 \SymColor{)}       }
\newcommand*{\WSRDVFB}[4][]{      \SymColor{\boldsymbol{H}}^{#1}         \SymColor{(} #2 , #3 ;#4\SymColor{)}       }

\newcommand*{\biasWSRDVmat}[2][]{   \SymColor{\boldsymbol{F}_{\SymColor{\mathrm{q}}}^{#1}}    \SymColor{(} #2 \SymColor{)}       }
\newcommand*{\biasWSRDFBgain}[1]{   \SymColor{\boldsymbol{G}_{\SymColor{\mathrm{q}}}}         \SymColor{(} #1 \SymColor{)}       }

\newcommand*{\ImpRSR}[2]{   \SymColor{\boldsymbol{f}}  \SymColor{(}  #1 \SymColor{,} #2 \SymColor{)}       }
\newcommand*{\ImpFBgain}[2]{   \SymColor{\boldsymbol{g}}  \SymColor{(}  #1 \SymColor{,} #2 \SymColor{)}       }
\newcommand*{\ImpTwoF}[2]{     \NonArgImpTwoF    \SymColor{(}  #1 \SymColor{,} #2 \SymColor{)}          }
\newcommand*{\NonArgImpTwoF}{   \SymColor{\boldsymbol{h}}     }

\newcommand*{\Vmat}{\SymColor{\boldsymbol{\Pi}}}
\newcommand*{\optVmat}[1]{       \Vmat_{#1}            }
\newcommand*{\optFBgain}[1]{     \FBgain_{#1}           }

\newcommand*{\IoptVmat}{     {\optVmat{\SymColor{\ast}}}            }
\newcommand*{\IoptFBgain}{   {\optFBgain{\SymColor{\ast}}}           }
\newcommand*{\IestVmat}{     \widehat{\optVmat{\SymColor{\ast}}}            }
\newcommand*{\IestFBgain}{   \widehat{\optFBgain{\SymColor{\ast}}}           }
\newcommand*{\IiniVmat}{     {\optVmat{0}}            }
\newcommand*{\IiniFBgain}{   {\optFBgain{0}}           }

\newcommand*{\KoIoptVmat}{     {\optVmat{\SymColor{\ast\ast}}}            }
\newcommand*{\KoIoptFBgain}{   {\optFBgain{\SymColor{\ast\ast}}}           }

\newcommand*{\SolPair}{   \SymColor{\boldsymbol{z}}       }
\newcommand*{\iteSolPair}[1]{   \SymColor{\boldsymbol{z}}_{#1}       } 
\newcommand*{\optSolPair}[1]{   \SymColor{\boldsymbol{z}}_{\SymColor{\ast}} \SymColor{(}  #1  \SymColor{)}       }


\newcommand*{\Ltrans}[4]{        \SymColor{\mathcal{T}}_{\SymColor{\mathrm{q}},#2}^{#4}         \SymColor{(} #3 \SymColor{)}       }
\newcommand*{\NonArgLtrans}[2]{  \SymColor{\mathcal{T}}_{\SymColor{\mathrm{q}},#1}^{#2}  }

\newcommand*{\kronLtrans}[1]{        \SymColor{\boldsymbol{M}}_{\SymColor{\mathrm{q}},#1}  }


\newcommand*{\smoothLBriskParam}{   \SymColor{\underline{\riskParam}_{\mathrm{s}}}    }
\newcommand*{\smoothUBriskParam}{   \SymColor{\overline{\riskParam}_{\mathrm{s}}}    }

\newcommand*{\uniqueLBriskParam}{   \SymColor{\underline{\riskParam}_{\mathrm{u}}}    }
\newcommand*{\uniqueUBriskParam}{   \SymColor{\overline{\riskParam}_{\mathrm{u}}}    }

\newcommand*{\NewtonLBriskParam}{   \SymColor{\underline{\riskParam}_{\mathrm{n}}}    }
\newcommand*{\NewtonUBriskParam}{   \SymColor{\overline{\riskParam}_{\mathrm{n}}}    }
\newcommand*{\NewtonLLBriskParam}{   \SymColor{\underline{\riskParam}_{\mathrm{n}}^{\prime}}    }
\newcommand*{\NewtonUUBriskParam}{   \SymColor{\overline{\riskParam}_{\mathrm{n}}^{\prime}}    }

\newcommand*{\StableLBriskParam}{   \SymColor{\underline{\riskParam}_{\mathrm{ms}}}    }
\newcommand*{\StableUBriskParam}{   \SymColor{\overline{\riskParam}_{\mathrm{ms}}}    }

\newcommand*{\NewtonConst}{\SymColor{K}_{\SymColor{\mathrm{N}}}}

\newcommand*{\MyPercentWorse}{\SymColor{\rho}}

\newcommand*{\SSOCPsetting}{\SymColor{ {|_{w=1}^{\mathrm{q}}} }}

\newcommand*{\MyEpsilonForConv}{\SymColor{\epsilon}}
\newcommand*{\MyDeltaForConv}{\SymColor{\delta}}


\newcommand*{\NewtonLipschitz}{   \SymColor{\gamma}_{\partial\NonArgImpTwoF{}{}}    }
\newcommand*{\NewtonScalar}{   \SymColor{\gamma^{\prime}}    }

\newcommand*{\DefBallRadius}{   \SymColor{\delta}    }

\newcommand*{\reqBallRadius}{   \SymColor{\delta_{\ast}}    }
\newcommand*{\LBreqBallRadius}{   \SymColor{\underline{\delta}_{\ast}}    }


\newcommand*{\COEFpsdMatCompare}{\SymColor{\kappa}}
\newcommand*{\unitVec}[1]{  \SymColor{\boldsymbol{e}}_{#1}  }

\newcommand*{\BoundImpTwoF}{\SymColor{\delta_{h}}}
\newcommand*{\BoundpartialImpTwoF}{\SymColor{\delta_{\partial h}}}
\newcommand*{\tmpUBriskParam}{\SymColor{\overline{\theta}}}
\newcommand*{\tmpLBriskParam}{\SymColor{\underline{\theta}}}

\newcommand*{\QuadFBMat}{\SymColor{\boldsymbol{\Psi}}}
\newcommand*{\EigQuadFBMat}[2]{\SymColor{\lambda_{#1}(#2)}}

\newcommand{\MyProof}{proof}

\newcommand{\MyFigRePos}[2]{\begin{center}\begin{minipage}[b]{#1} #2 \end{minipage}\end{center}}

\section{Introduction} \label{sec_intro}

Noise and uncertainty contained in dynamical systems are expressed by stochastic system parameters \citep{Mesbah16}.
Independent and identically distributed (i.i.d.) stochastic parameters such as those in \citep{Koning82} have attracted significant attention because they can represent various noises and uncertainties.
For example, the practical applications of the i.i.d. parameters involve sensorimotor systems \citep{Todorov05}, time-varying communication delays in networks \citep{Hosoe22CDC}, vehicle platoons via lossy communication \citep{Acciani22}, and digital control with random sampling intervals \citep{Koning82}.
Input- and state-dependent noise in aerospace systems \citep{Mclane71} can be represented using i.i.d. parameters with the discretization of the systems. 	
I.i.d. parameters have been extended to combinations with other stochastic parameters \citep{fujisaki07,Fisher09,ItoCybern23} to treat complex uncertainties \citep{HosoeAutoma20}.

Several stability notions and optimal control policies have been proposed for linear systems with i.i.d. stochastic parameters.
Stochastic optimal control \citep{Koning82} minimizes an average of cost functions by extending traditional deterministic optimal control \citep{Anderson89}.
Another stochastic optimal control has been developed to suppress the variance of system states \citep{Fujimoto11c}.
These control laws guarantee mean and mean-square (MS) stabilities, which are asymptotic stabilities in the first- and second-order moments, respectively.
Stability of high-order moments has been analyzed \citep{Luo20TAC,Zhang2020AMC,Zhang2021SCIS,Zhang2022TSMC,Ogura13}, which includes both notions of stability and robustness to system randomness \citep{ItoCSL23}.

A crucial challenge is to design multiple types of controllers to handle statistical properties of systems with i.i.d. stochastic parameters.
However, such complex controllers are difficult to design appropriately.
Risk-sensitive (RS) control is a promising example to handle the risk of unexpected control results \citep{Jacobson73,Duncan13,Lim04,Medina12}.
Even if the aforementioned stochastic optimal control realizes the desired average performance, worse-case results are often critical and should be avoided.
RS control with a risk-averse policy is helpful in mitigating worse results rather than average results.
By contrast, RS control with a risk-seeking policy specializes in enhancing better results.
However, designing RS controllers for linear systems with i.i.d. parameters is not straightforward.
These controllers can be nonlinear, whereas linear controllers are highly compatible with linear systems in terms of reliability and implementation.
Because such nonlinearity makes the controller design difficult, the design has relied on approximation methods \citep{Ruszczynski10,Shen14,Broek10}.
A trade-off exists between the size of the state region and risk sensitivity \citep{Nagai95}.
The details are discussed in \citep{ItoFHRSL19}.

An underlying difficulty in linear i.i.d. systems is that Riccati-like algebraic equations cannot be employed to design complex controllers.
Although classical RS controllers are designed based on algebraic equations \citep{Jacobson73}, they are violated if i.i.d. parameters are included.
Our previous work \citep{ItoFHRSL19} has addressed this difficulty and proposed risk-sensitive linear (RSL) control for linear i.i.d. systems.
The RSL control overcomes the aforementioned drawbacks and realizes the following: the controller is linear; its exact solution is derived; and it operates on the entire state space. 
Nonetheless, the following problems remain.
Our previous work has focused on RS control without considering the possibility of designing more general types of controllers. 
The design of RSL controllers over an infinite-horizon (IH) case remains challenging while its concept has been presented.
Specifically, 
the IH-RSL controller design incurs a huge computational cost via the iteration of solving nonlinear optimization, 
and
stability and optimality of the IH-RSL control should be theoretically guaranteed.

To overcome the aforementioned problems, this paper presents a general framework for designing various linear optimal controllers, including RSL controllers.
We propose weighted stochastic Riccati (WSR) equations, which are powerful tools for designing IH controllers for linear systems with i.i.d. stochastic parameters.
The main contributions of this study are summarized as follows.

\begin{enumerate}
	
	\item	
	\textbf{Generality:} 
	Solving the proposed WSR equations is shown to provide multiple types of IH optimal controllers by designing a weight function (Theorem \ref{thm:solution_problem_1}).
	The existing stochastic optimal \citep{Koning82}, RSL \citep{ItoFHRSL19}, and novel RS controllers are covered (Examples \ref{rem:SSOC}--\ref{rem:RRSL} in Section \ref{sec_main_problems}).
	
	\item
	\textbf{Solvability:} 
	We propose two approaches to solve the WSR equations.
	The first approach is to derive WSR difference equations.
	Iterative solutions to the WSR difference equations converge to solutions to the WSR equations (Theorem \ref{thm:solution_WSR}).
	In the second approach, Newton's method is employed.
	We derive proper initialization needed for using Newton's method (Theorem \ref{thm:MyNewton}). 
	Moreover, 
	we show the uniqueness and smoothness of the solution to the WSR equations (Theorem \ref{thm:uniqueness}).

	\item
	\textbf{Novel control:}
	As one example of using the proposed framework, we propose robust RSL (RRSL) controllers (Example \ref{rem:RRSL} in Section \ref{sec_main_problems}).
	While the RRSL controllers enable the realization of an RS control policy,
	they are more robust than the existing RSL controllers in terms of the randomness occurring in the stochastic controller design.
	In other words, the controller design often needs to approximate expectations regarding i.i.d. parameters by using random samples.
	The RRSL controllers suppress the degradation of the design caused by such random samples.

	\item
	\textbf{Advantages:}
	The proposed design using the WSR equations has the following advantages compared with the previous design \citep{ItoFHRSL19}.
	Stability of the feedback system with applying the designed controller is guaranteed (Section \ref{sec_method2}).
	Optimality is guaranteed for IH controllers rather than finite-horizon controllers (Theorem \ref{thm:solution_problem_1}).
	The computational cost of the proposed design is reduced because it does not need iterative nonlinear optimization whereas the previous design needs it (Section \ref{sec_solve_problem_1_solve}).

	\item
	\textbf{Demonstration:}
	Numerical examples are presented to show the effectiveness of the proposed method in terms of convergence, robustness, stability, and control performance (Section \ref{sec_numerical_example}).
	
\end{enumerate}

This paper is a substantially extended version of our conference paper \citep{ItoACC16} and its main extensions are summarized below.
This study proposes the WSR equations associated with theoretical analyses, which are generalized versions of limited equations in the conference paper.
An analysis of the WSR equations and approach based on Newton's method are additionally presented for solving the WSR equations.
This study proposes novel RRSL controllers that are more robust than RSL controllers presented in the conference paper.	
Several stability analyses are additionally presented.
All numerical simulations are novel materials to demonstrate the effectiveness of the proposed method.

The remainder of this paper is organized as follows.
Section \ref{sec_problem_setting} gives two main problems in this study.
Our solutions to these two problems are proposed in Sections \ref{sec_method} and \ref{sec_method2}. 
In Section \ref{sec_numerical_example}, the effectiveness of the proposed method is evaluated using numerical simulations.
Finally, this study is concluded in Section \ref{sec_conclusion}.

\textbf{Notation: }
The following notations are used:
\begin{itemize}

	\item
	${\SetSymMat{\DimANotation}}$: the set of $\DimANotation \times \DimANotation$ real-valued symmetric matrices

	\item 
	$\Identity{\DimANotation}$: the $\DimANotation \times \DimANotation$ identity matrix
	\item
	$\El{\NotationVec}{\IDNotation}$: the $\IDNotation$-th component of a vector $\NotationVec \in \mathbb{R}^{\DimANotation}$
	\item
	$\El{\NotationMat}{\IDNotation,\IDbNotation}$: the component in the $\IDNotation$-th row and $\IDbNotation$-th column  of a matrix $\NotationMat \in \mathbb{R}^{\DimANotation \times \DimBNotation}$

	\item
	$\VEC{\NotationMat}:=[ 
	\El{\NotationMat}{1,1} , \dots, \El{\NotationMat}{\DimANotation,1} , 
	\El{\NotationMat}{1,2} , \dots, \El{\NotationMat}{\DimANotation,2} , 
	\dots ,
	$ $	
	\El{\NotationMat}{1,\DimBNotation} , \dots, \El{\NotationMat}{\DimANotation,\DimBNotation}	]^{\MyTRANSPO} $: 
	the vectorization of the components of a matrix $\NotationMat \in \mathbb{R}^{\DimANotation \times \DimBNotation}$
	
	\item
	$\VECH{\NotationSquMat}:=[ 
	\El{\NotationSquMat}{1,1} , \dots, \El{\NotationSquMat}{\DimANotation,1} , 
	\El{\NotationSquMat}{2,2} , \dots, \El{\NotationSquMat}{\DimANotation,2} , 
	\dots ,
	$ $
	\El{\NotationSquMat}{\IDbNotation,\IDbNotation} , \dots, \El{\NotationSquMat}{\DimANotation,\IDbNotation} , 
	\dots ,
	\El{\NotationSquMat}{\DimANotation,\DimANotation}	]^{\MyTRANSPO} $: 
	the half vectorization of the lower triangular components of a square matrix $\NotationSquMat \in \mathbb{R}^{\DimANotation \times \DimANotation}$

	\item
	$\NotationMat_{\NotationA} \otimes \NotationMat_{\NotationB} \in \mathbb{R}^{\DimANotation_{\NotationA} \DimANotation_{\NotationB} \times \DimBNotation_{\NotationA} \DimBNotation_{\NotationB} }$: the Kronecker product of matrices $\NotationMat_{\NotationA} \in \mathbb{R}^{\DimANotation_{\NotationA} \times \DimBNotation_{\NotationA} }$ and $\NotationMat_{\NotationB} \in \mathbb{R}^{\DimANotation_{\NotationB} \times \DimBNotation_{\NotationB} }$, given by
	\begin{align*}
		\NotationMat_{\NotationA} \otimes \NotationMat_{\NotationB} =
		\begin{bmatrix}
			\El{\NotationMat_{\NotationA}}{1,1} \NotationMat_{\NotationB} & \hdots & \El{\NotationMat_{\NotationA}}{1,\DimBNotation_{\NotationA}} \NotationMat_{\NotationB} \\
			\vdots & \ddots & \vdots \\
			\El{\NotationMat_{\NotationA}}{\DimANotation_{\NotationA},1} \NotationMat_{\NotationB} & \hdots & \El{\NotationMat_{\NotationA}}{\DimANotation_{\NotationA},\DimBNotation_{\NotationA}} \NotationMat_{\NotationB} 
		\end{bmatrix}
		,			
	\end{align*}	
	where $\NotationMat^{\otimes 2}:= \NotationMat \otimes \NotationMat$

	\item
	${\DupMat{\DimANotation}}$, ${\EliMat{\DimANotation}}$: the duplication matrix and elimination matrix that satisfy 
	${\DupMat{\DimANotation}} \VECH{\NotationSymMat}=\VEC{\NotationSymMat}$, ${\EliMat{\DimANotation}} \VEC{\NotationSymMat}=\VECH{\NotationSymMat}$, and ${\EliMat{\DimANotation}}{\DupMat{\DimANotation}}=\Identity{\DimANotation(\DimANotation+1)/2}$ for any symmetric matrix $\NotationSymMat \in  {\SetSymMat{\DimANotation}}$ \citep[Definitions 3.1a, 3.1b, 3.2a, and 3.2b and Lemma 3.5 (i)]{Magnus80}.
	The examples for $\DimANotation=2$ are as follows:
	\begin{align*}		
		{\EliMat{2}}
		&=
		\begin{bmatrix}
			1 & 0 & 0 & 0 \\
			0 & 1 & 0 & 0 \\
			0 & 0 & 0 & 1 
		\end{bmatrix}
		,\quad
		{\DupMat{2}}
		=
		\begin{bmatrix}
			1 & 0 & 0 \\
			0 & 1 & 0 \\
			0 & 1 & 0 \\
			0 & 0 & 1 
		\end{bmatrix}
		.
	\end{align*}

	\item
	$\NotationSymMat \succ 0$ (resp. $\prec 0$): the positive (resp. negative) definiteness of a symmetric
	matrix $\NotationSymMat \in {\SetSymMat{\DimANotation}}$
	\item
	$\NotationSymMat \succeq 0$ (resp. $\preceq 0$): the positive (resp. negative) semidefiniteness of a symmetric matrix $\NotationSymMat \in {\SetSymMat{\DimANotation}}$

	\item 
	$\NotationVec_{\IDEl:\IDbEl}:=[\NotationVec_{\IDEl} , \NotationVec_{\IDEl+1} \dots , \NotationVec_{\IDbEl}]$: the array of $\NotationVec_{\MyT}$ for $\IDEl \leq \MyT \leq \IDbEl$
	
	\item
	$\partial_{\NotationVec^{\MyTRANSPO}} \NotationFunc (\NotationVec)={\partial \NotationFunc (\NotationVec)}/{\partial \NotationVec^{\MyTRANSPO}}$: the partial derivative of $\NotationFunc(\NotationVec)$ with respect to $\NotationVec$, indicating
	${\El{ {\partial \NotationFunc (\NotationVec)}/{\partial \NotationVec^{\MyTRANSPO}} }{\IDEl,\IDbEl}}
	={\partial {\El{ \NotationFunc(\NotationVec)}{\IDEl}} }/
	{\partial {\El{ \NotationVec}{\IDbEl}} }$

	\item 
	${\Expect{\StoParam}}[\NotationFunc(\StoParam)]$: the expectation of a function $\NotationFunc(\StoParam)$ with respect to a random vector $\StoParam$
	
	\item 
	$\Covariance{\StoParam}[\NotationFunc(\StoParam)]$: the covariance of a vector-valued function $\NotationFunc(\StoParam)$ with respect to a random vector $\StoParam$

\end{itemize}

\section{Problem setting} \label{sec_problem_setting}

Section \ref{sec_target_sys} describes the target systems considered in this study.
Main problems to these systems are presented in Section \ref{sec_main_problems}.

\subsection{Target systems}\label{sec_target_sys}

Let us consider the linear system with stochastic system matrices:
\begin{align}
	\State_{\MyT+1}&\;=\DriftMat_{\MyT} \State_{\MyT} + \InMat_{\MyT} \Input_{\MyT}
	,\label{eq:def_sys}
	\\ 
	\StoParam_{\MyT}&:=[ \VEC{\DriftMat_{\MyT}}^{\MyTRANSPO}, \VEC{\InMat_{\MyT}}^{\MyTRANSPO}]^{\MyTRANSPO}
	,\label{eq:def_Lambda}
\end{align}	
where $\State_{\MyT} \in \mathbb{R}^{\DimX}$ and $\Input_{\MyT} \in \mathbb{R}^{\DimU}$ are the state and control input at the time $\MyT$, respectively.
The initial state $\State_{0}$ is deterministic.
The stochastic parameter $\StoParam_{\MyT} \in \DomStoParam \subseteq \mathbb{R}^{\DimX(\DimX+\DimU)}$ denotes the array of the stochastic system matrices $\DriftMat_{\MyT} \in \mathbb{R}^{\DimX \times \DimX}$ and $\InMat_{\MyT} \in \mathbb{R}^{\DimX \times \DimU}$.
Assume that the probability density function (PDF)
$\PDF{\StoParam}$ of $\StoParam$ is known and that $\StoParam_{\MyT}$ is an i.i.d. random sample of $\StoParam$ at each time $\MyT$.
The PDF $\PDF{\StoParam}$ is assumed to be a continuous function on a Lebesgue measurable set $\DomStoParam$.

\subsection{Main problems}\label{sec_main_problems}

We consider the following feedback system by applying a state feedback controller $\funcInput: \mathbb{R}^{\DimX}\to \mathbb{R}^{\DimU}$ to \eqref{eq:def_sys}:
\begin{align}
	\State_{\MyT+1}&=\DriftMat_{\MyT} \State_{\MyT} + \InMat_{\MyT} \funcInput(\State_{\MyT})
	.\label{eq:def_FBsys}
\end{align}	
Introducing a weight function ${\Weight{}}$ associated with a sensitivity parameter $\riskParam \in \mathbb{R}$, we consider an IH version of the weighted cost function \citep{ItoFHRSL19}:
\begin{align}	
	{\IHWcostJ{\funcInput}{\State_{0}}{\riskParam}}
	&:= 
	\lim_{\terminalT \rightarrow \infty}
	{\Expect{\StoParam_{0:\terminalT}}}\Big[
	\sum_{\MyT = 0}^{\terminalT } 
	\Big(\Big(  \prod_{\MybT = 0}^{\MyT}  {\Weight{\MybT}} \Big) 	
	\nonumber\\&
	\qquad\times
	\Big( 
	\State_{\MyT}^{\MyTRANSPO} \xCostMat \State_{\MyT} 
	+ \funcInput(\State_{\MyT})^{\MyTRANSPO} \uCostMat \funcInput(\State_{\MyT}) 
	\Big)\Big) 
	\Big]
	,\label{eq:def_IHWcostJ}
\end{align}	
where $\uCostMat \succ 0 \in {\SetSymMat{\DimU}}$ and $\xCostMat \succ 0\in {\SetSymMat{\DimX}}$  are given positive definite matrices.
Various performance metrics are expressed according to the setting of ${\Weight{\MybT}}$.

Let us propose a desired weight function $\dWeight(\StoParam; \riskParam, \funcInput  ,  \funccostJ  )$ as a reference for ${\Weight{}}$.
The desired weight $\dWeight(\StoParam; \riskParam, \funcInput  ,  \funccostJ  )$ is a function of $(\StoParam, \riskParam)$ and is a functional of $(\funcInput, \funccostJ)$,
where $\funccostJ: \mathbb{R}^{\DimX}\to \mathbb{R}$ is an estimate  of the cost function ${\IHWcostJ{\funcInput}{\wildcard}{\riskParam}}$ given $\funcInput$ and $\riskParam$.
Throughout this study, we use the following assumption:
\begin{assumption}[Desired weight function]\label{ass:desired_weight}
	Given a sensitivity parameter $\riskParam$ and functions $\funcInput$ and $\funccostJ$, 
	a desired weight $\dWeight$ satisfies the following conditions:
	\begin{enumerate}
	\item 
		The desired weight $\dWeight(\StoParam; \riskParam, \funcInput  ,  \funccostJ  )$ is continuous in $\StoParam$ on $\DomStoParam$.	
		
		\item
		We have ${\Expect{\StoParam}}[ \dWeight(\StoParam; \riskParam, \funcInput  ,  \funccostJ  ) ]=1$
		and $ \dWeight(\StoParam; \riskParam, \funcInput  ,  \funccostJ  ) \geq 0$ for any $\StoParam \in \DomStoParam$.
		\item
		If $\riskParam=0$ holds, we have $\dWeight(\StoParam; \riskParam, \funcInput  ,  \funccostJ  )=1$ for any $\StoParam \in \DomStoParam$.

	\end{enumerate}
	
\end{assumption}

This study addresses two problems.
The first problem is as follows:

\textbf{Problem 1 (Controller design):}
Given a desired weight $\dWeight$ and sensitivity parameter $\riskParam$, find an optimal feedback controller $\optFBcontroller: \mathbb{R}^{\DimX}\to \mathbb{R}^{\DimU}$, minimum cost $\mincostJ(\State)$, and weight ${\Weight{}}$ that satisfy 
\begin{align}
	\forall \State_{0} \in \mathbb{R}^{\DimX},\;	
	\optFBcontroller
	&\in
	\argmin_{\funcInput} 
	{\IHWcostJ{\funcInput}{\State_{0}}{\riskParam}}
	, \label{eq:def_optFBcontroller}
	\\
	\forall \State_{0} \in \mathbb{R}^{\DimX},\;
	\mincostJ(\State_{0})
	&:=
	\min_{\funcInput} 
	{\IHWcostJ{\funcInput}{\State_{0}}{\riskParam}}
	,\label{eq:def_mincostJ}
	\\
	\forall \StoParam \in \DomStoParam,\;
	{\Weight{}}
	&=
	\dWeight(\StoParam; \riskParam, \optFBcontroller ,  \mincostJ )
	.\label{eq:problem1_Weight}
\end{align}

Various control policies can be considered in Problem 1 according to the setting of the desired weight $\dWeight$, examples of which are introduced below.

\begin{example}[Standard stochastic optimal control]\label{rem:SSOC}
	If we set
	$\dWeight=1$,
	we have ${\Weight{}}=1$ from \eqref{eq:problem1_Weight} and
	Problem 1 reduces to an existing stochastic optimal control problem involving \eqref{eq:def_optFBcontroller} and \eqref{eq:def_mincostJ} \citep{Koning82}.
\end{example}

\begin{example}[Risk-sensitive linear control]\label{rem:RSL}
	If we set $\dWeight$ as follows:
	\begin{align}	
		\dWeight(\StoParam; \riskParam, \funcInput  ,  \funccostJ  )
		&\propto
		\exp  \riskParam \predcostJ(\StoParam;  \funcInput  ,  \funccostJ  )
		,
		\label{eq:def_dWeight_RSL}	
		\\
		\predcostJ(\StoParam;  \funcInput  ,  \funccostJ  )
		&:=
		{\Expect{\RefX}}[
		\funccostJ(   \DriftMat \RefX + \InMat \funcInput(\RefX)  )
		\nonumber\\&\quad
		+	\RefX^{\MyTRANSPO} \xCostMat \RefX
		+ \funcInput(\RefX)^{\MyTRANSPO} \uCostMat \funcInput(\RefX)  
		],
		\label{eq:def_predcostJ}	
	\end{align}
	then Problem 1 can be interpreted as a {slightly modified} version of the IH RSL control problem \citep{ItoFHRSL19}.
	Solving Problem 1 with the desired weight \eqref{eq:def_dWeight_RSL} yields an RS controller.
	The control policy depends on $\riskParam$; setting $\riskParam>0$ leads to risk-averse control to mitigate worse cases in various control results.
	For a controller $\funcInput$, $\predcostJ(\StoParam;  \funcInput  ,  \funccostJ  )$
	indicates the predictive cost of the per-step state transition $\DriftMat \RefX + \InMat \funcInput(\RefX)$ expected over a random state $\RefX \in \mathbb{R}^{\DimX}$.	
	In \eqref{eq:def_dWeight_RSL}, the exponential form of $\predcostJ(\StoParam;  \funcInput  ,  \funccostJ  )$ acts as a risk measure for the per-step state transition depending on each $\StoParam$.
	Further analyses of the RSL control are presented in \citep{ItoFHRSL19}.
\end{example}

\begin{example}[Robust risk-sensitive linear control]\label{rem:RRSL}
	We propose RRSL controllers to enhance robustness to randomness occurring in the RSL controller design.
	For general PDFs of $\StoParam$, the expectations ${\Expect{\StoParam_{\MyT}}}[\dots]$ included in equations used for the design are often approximated using the Monte Calro (MC) method with random samples of $\StoParam$.
	The RRSL controllers employ the following desired weight, which is robust to sample randomness: 
	\begin{align}	
		&\dWeight(\StoParam; \riskParam, \funcInput  ,  \funccostJ  )
		\nonumber\\&
		\propto
		1
		+	
		\frac{\riskParam}{1 + \exp( 
			- \paramWeightA   \predcostJ(\StoParam; \funcInput  ,  \funccostJ  ) 
			+ \paramWeightB {\Expect{\StoParam}}[ \predcostJ(\StoParam; \funcInput  ,  \funccostJ  )]
			)  }
		,\label{eq:def_dWeight_RRSL}	
	\end{align}
	where $\paramWeightA \in \mathbb{R}$ and $\paramWeightB \in \mathbb{R}$ are free parameters. 
	While this weight enables the realization of an RS control policy, it is more robust than the weight \eqref{eq:def_dWeight_RSL} of the RSL control with $\riskParam>0$.
	Intuitively, we obtain the weight \eqref{eq:def_dWeight_RRSL} by replacing the exponential function in \eqref{eq:def_dWeight_RSL} with the sigmoid function.
	As illustrated in Fig. \ref{fig:exp_vs_sig} \subref{fig:exp_vs_sig_func},
	for a one-dimensional $\StoParam$ and desired function $\dWeight(\StoParam)$,
	both sigmoid and exponential $\dWeight$ values emphasize higher values of $\StoParam$.
	The histograms in Fig. \ref{fig:exp_vs_sig} \subref{fig:exp_vs_sig_rnd} show that the sigmoid values of $\StoParam$ are hardly dispersed in comparison with the exponential values.
	This indicates the robustness of the RRSL control for random samples.
		The RRSL control focuses on a risk-averse case ($\riskParam>0$) while the randomness is not serious in a risk-seeking case ($\riskParam<0$) in which the exponential function in \eqref{eq:def_dWeight_RSL} does not diverge.
	The effectiveness of the RRSL control is demonstrated in Section \ref{sec_numerical_example}. 
\end{example}

\definecolor{mycolorExpBin}{rgb}{ 0.4,0.4,0.52 }
\definecolor{mycolorSigBin}{rgb}{ 1.0,0.64,0.64 }

\begin{figure}[t]\MyFigRePos{.55\linewidth}{	
	\begin{minipage}[b]{.5\linewidth}
		\centering	
		{\footnotesize
			\hfill
			\begin{tabular}{| l |}
				\hline
				{\textcolor{blue}{-\;-}}: Exponential $\dWeight$
				\\
				{\textcolor{red}{---}}: Sigmoid $\dWeight$
				\\
				\hline			
			\end{tabular}	
		}
		\begin{tikzpicture}[scale=1.0]
		\begin{axis}[axis y line=none, axis x line=none
		,xtick=\empty, ytick=\empty
		,xmin=-10,xmax=10,ymin=-10,ymax=10
		,width=2.25 in,height=1.55 in
		]
		
		\node at (0,-0.2) { 
			\includegraphics[width=0.95\linewidth]{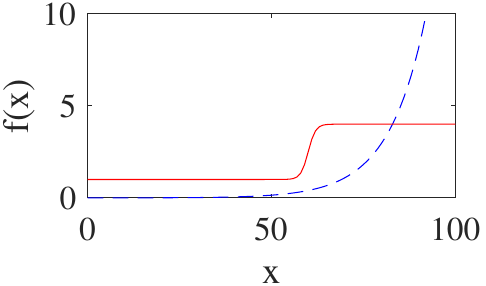}
		};

		\filldraw[draw=white, fill=white] (-15,-15) -- (-15,15) -- (-8,15) -- (-8,-15);	
		\filldraw[draw=white, fill=white] (-15,-15) -- (-15,-7.6) -- (15,-7.6) -- (15,-15);	
		
		\node at (-8.9,1) {\rotatebox{90}{\small	
				{\shortstack{Value of $\dWeight(\StoParam)$}} 
		}};
		\node at (1,-8.7) {\rotatebox{0}{\small	
				{\shortstack{Value of $\StoParam$}} 
		}};
		
		\end{axis}			
		\end{tikzpicture}		
		\subcaption{An exponential and a sigmoid function.}\label{fig:exp_vs_sig_func}
	\end{minipage}%
	\begin{minipage}[b]{.5\linewidth}
		\centering	
		{\footnotesize
			\hfill
			\begin{tabular}{| l |}
				\hline
				{\textcolor{mycolorExpBin}{$\blacksquare$}}: Exponential $\dWeight$
				\\
				{\textcolor{mycolorSigBin}{$\blacksquare$}}: Sigmoid $\dWeight$
				\\
				\hline			
			\end{tabular}	
		}
		\begin{tikzpicture}[scale=1.0]
		\begin{axis}[axis y line=none, axis x line=none
		,xtick=\empty, ytick=\empty
		,xmin=-10,xmax=10,ymin=-10,ymax=10
		,width=2.25 in,height=1.55 in
		]
	
		\node at (0,-0.2) { 
			\includegraphics[width=0.92\linewidth]{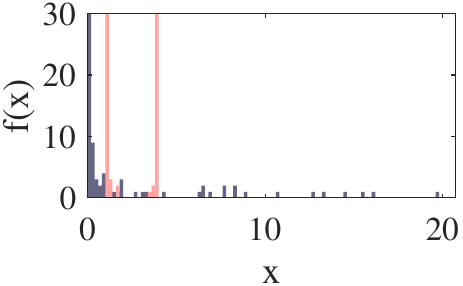}
		};

		\filldraw[draw=white, fill=white] (-15,-15) -- (-15,15) -- (-8,15) -- (-8,-15);	
		\filldraw[draw=white, fill=white] (-15,-15) -- (-15,-7.6) -- (15,-7.6) -- (15,-15);	
		\node at (-8.8,2) {\rotatebox{90}{\small	
				{\shortstack{Frequency}} 
		}};
		\node at (1,-8.7) {\rotatebox{0}{\small	
				{\shortstack{Value of $\dWeight(\StoParam)$}} 
		}};
		
		\end{axis}			
		\end{tikzpicture}		
		\subcaption{Histograms generated by randomly sampled $\StoParam$}\label{fig:exp_vs_sig_rnd}
	\end{minipage}%
	\caption{Comparison between exponential and sigmoid functions $\dWeight(\StoParam)$ from the perspective of robustness to random samples $\StoParam$.}\label{fig:exp_vs_sig}
}\end{figure}
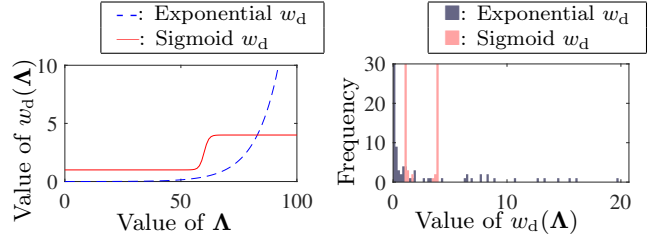

The second problem focuses on the well-known MS stability of the system with controllers obtained by solving Problem 1.
This study proposes the following generalized version of the MS stability.
\begin{definition}[Weighted mean square stability]\label{def:wms_stable}
	Given weight functions ${\WmsWeight{\MyT}}  \in \mathbb{R}$ for $\MyT \in \{0,1,2,\dots\}$,
	the feedback system \eqref{eq:def_FBsys} is said to be \textit{weighted mean square (WMS) stable} with ${\WmsWeight{\MyT}}$ 
	if for each $\State_{0} \in \mathbb{R}^{\DimX}$, we have $\lim_{\MyT \rightarrow \infty} {\Expect{\StoParam_{0:\MyT}}}[  {\WmsWeight{\MyT}}   \|\State_{\MyT}\|^{2} ]=0$.
	The WMS stability with ${\WmsWeight{\MyT}}=1$ is equivalent to the MS stability \citep{Koning82}.
	The feedback system \eqref{eq:def_FBsys} is said to be MS (resp. WMS) stabilizable if there exists $\FBgain \in \mathbb{R}^{\DimU \times \DimX}$ such that the system with $\funcInput(\State)=-\FBgain \State$ is MS (resp. WMS) stable.
\end{definition}

\textbf{Problem 2 (Second-moment stability analysis):} 
Analyze the WMS stability of the feedback system \eqref{eq:def_FBsys} with applying $\optFBcontroller$ that is obtained by solving Problem 1.

\section{Proposed method: solution to Problem 1}\label{sec_method}

Solving Problem 1 reduces to solving the WSR equations proposed in Section \ref{sec_solve_problem_1}. 
Section \ref{sec_solve_problem_1_solve} presents two approaches for solving the WSR equations.

\subsection{Derivation of the WSR equations}\label{sec_solve_problem_1}

We propose the WSR equations that are key to this study below.

\begin{definition}[WSR equations]
	Given a desired weight $\dWeight$ and sensitivity parameter $\riskParam$,
	let us define the WSR equations of $({\Vmat},{\FBgain}) \in \SetSymMat{\DimX} \times \mathbb{R}^{\DimU \times \DimX} $:
	\begin{align}	
		\Vmat
		&={\WSRDVmat{\Vmat}{\FBgain}{\riskParam}}
		,\label{eq:def_optVmat_inf}
		\\
		\FBgain
		&={\WSRDFBgain{\Vmat}{\FBgain}{\riskParam}}	
		,\label{eq:def_optFBgain_inf}
	\end{align}		
	where
	${\WSRDVmat{\wildcard}{\wildcard}{\riskParam}} :\DomALLVmat \times \DomALLFBgain \to 	\SetSymMat{\DimX}$
	and
	${\WSRDFBgain{\wildcard}{\wildcard}{\riskParam}} :\DomALLVmat \times \DomALLFBgain\to  \mathbb{R}^{\DimU \times \DimX}$
	are given by	
	\begin{align} 
		{\WSRDVmat{\Vmat}{\FBgain}{\riskParam}}		
		&:=
		{\WE{\StoParam}{\riskParam}{\FBgain}{\Vmat}}
		[ \DriftMat^{\MyTRANSPO} \Vmat \DriftMat] 
		+
		\xCostMat
		\nonumber\\&
		\quad
		-
		{\WE{\StoParam}{\riskParam}{\FBgain}{\Vmat}}
		[ \DriftMat^{\MyTRANSPO} \Vmat \InMat] 
		{\WSRDFBgain{\Vmat}{\FBgain}{\riskParam}}			
		,\label{eq:def_WSRDVmat}
		\\
		{\WSRDFBgain{\Vmat}{\FBgain}{\riskParam}}	
		&
		:=
		{\WE{\StoParam}{\riskParam}{\FBgain}{\Vmat}}
		[ \InMat^{\MyTRANSPO} \Vmat \InMat  +  \uCostMat]^{-1}
		{\WE{\StoParam}{\riskParam}{\FBgain}{\Vmat}}
		[ \InMat^{\MyTRANSPO} \Vmat \DriftMat] 
		, \label{eq:def_WSRDFBgain}
	\end{align}
	and ${\WE{\StoParam}{\riskParam}{\FBgain}{\Vmat}}[\dots]$ denotes
	the weighted expectation: for any continuous function $\NotationFunc: \DomStoParam \to \mathbb{R}$,
	\begin{align}	
		{\WE{\StoParam}{\riskParam}{\FBgain}{\Vmat}}
		[ \NotationFunc(\StoParam) ]
		&:=
		{\Expect{\StoParam}}
		\big[ 
		\NotationFunc(\StoParam)
		\dWeight(\StoParam; \riskParam, \estfuncInput(\wildcard;\FBgain)  ,  \estfunccostJ(\wildcard;\Vmat)  )
		\big]
		,\label{eq:def_WE}
		\\
		\estfuncInput(\State;\FBgain)
		&:=  -\FBgain \State
		,\\
		\estfunccostJ(\State;\FBgain)
		&:=\State^{\MyTRANSPO} \Vmat \State
		,
	\end{align}	
	where $\DomALLVmat \times \DomALLFBgain$ is a subset of $\SetSymMat{\DimX} \times \mathbb{R}^{\DimU \times \DimX}$ on which  $\|{\WE{\StoParam}{\riskParam}{\FBgain}{\Vmat}}	[ \StoParam \StoParam^{\MyTRANSPO}]\| < \infty$ and ${\WE{\StoParam}{\riskParam}{\FBgain}{\Vmat}}
	[ \InMat^{\MyTRANSPO} \Vmat \InMat  +  \uCostMat]\succ 0$ are satisfied so that \eqref{eq:def_WSRDVmat} and \eqref{eq:def_WSRDFBgain} are well defined.

\end{definition}

\begin{remark}[Special cases]\label{rem:special_cases}
	If $\riskParam=0$, that is, $\dWeight=1$ holds, the WSR equations reduce to a stochastic version of the discrete-time Riccati equations \citep{Koning82}, as introduced in Appendix \ref{sec_supporting}.
	If $\dWeight=1$ holds and $(\DriftMat,\InMat)$ is deterministic, then the WSR equations are equivalent to the well-known deterministic Riccati equations. 
\end{remark}

\begin{theorem}[Solution to Problem 1]\label{thm:solution_problem_1}
	For any solution $(\IoptVmat,\IoptFBgain)$ to the WSR equations \eqref{eq:def_optVmat_inf} and \eqref{eq:def_optFBgain_inf} that satisfies $\IoptVmat \succ 0$,
	a solution to Problem 1 is given by
	\begin{align}	
		\optFBcontroller(\State;\IoptFBgain)&=-\IoptFBgain \State
		,\label{eq:sol_optFBcontroller}
		\\
		\mincostJ(\State;\IoptVmat)&= \State^{\MyTRANSPO} \IoptVmat \State 
		,\label{eq:sol_ValF}
	\end{align}
	where
	${\Weight{}}
	=
	\dWeight(\StoParam; \riskParam, \optFBcontroller(\wildcard;\IoptFBgain) , $ $  \mincostJ(\wildcard;\IoptVmat) )$
	in 	\eqref{eq:problem1_Weight}.
\end{theorem}
\begin{\MyProof}	
	The proof is described in Appendix \ref{pf:solution_problem_1}.	
\end{\MyProof}	

\begin{remark}[Contribution of Theorem \ref{thm:solution_problem_1}]
	Theorem \ref{thm:solution_problem_1} indicates that solving Problem 1 reduces to solving the WSR equations.
	Although the WSR equations are not algebraic in $(\Vmat,\FBgain)$ except for special cases, 
	we propose methods to solve them in Section \ref{sec_solve_problem_1_solve}.
\end{remark}

Next, we analyze the uniqueness and smoothness of a solution to the WSR equations.
Let $\SolPair$ be the vectorization of $(\Vmat,\FBgain)$,
and we consider the implicit form ${\ImpTwoF{\SolPair}{\riskParam}}=0$ of the WSR equations \eqref{eq:def_optVmat_inf} and \eqref{eq:def_optFBgain_inf} as follows:
\begin{align}
	\SolPair&:=
	\begin{bmatrix}
		\VECH{\Vmat}
		\\ 
		\VEC{\FBgain} 
	\end{bmatrix}
	\in \mathbb{R}^{ (\DimX(\DimX+1)/2) + \DimU \DimX   }
	,\\
	{\ImpTwoF{\SolPair}{\riskParam}}
	&:=
	\begin{bmatrix}
		{\ImpRSR{\SolPair}{\riskParam}}
		\\
		{\ImpFBgain{\SolPair}{\riskParam}}
	\end{bmatrix}
	\in \mathbb{R}^{ (\DimX(\DimX+1)/2) + \DimU \DimX   }
	,\\
	{\ImpRSR{\SolPair}{\riskParam}}
	&:=
	\VECH[\big]{
		{\WE{\StoParam}{\riskParam}{\FBgain}{\Vmat}}
		[{\DescriptionFBMat{\FBgain}}^{\MyTRANSPO} \Vmat {\DescriptionFBMat{\FBgain}}]
		\nonumber\\&\quad
		+  \FBgain^{\MyTRANSPO} \uCostMat \FBgain   
		+ \xCostMat
		-\Vmat	
	}
	\in \mathbb{R}^{ \DimX(\DimX+1)/2 }
	,\\
	{\ImpFBgain{\SolPair}{\riskParam}}
	&:=
	\VEC[\big]{
		{\WE{\StoParam}{\riskParam}{\FBgain}{\Vmat}}
		[ \InMat^{\MyTRANSPO} \Vmat \InMat  + \uCostMat]\FBgain
		\nonumber\\&\quad
		-	
		{\WE{\StoParam}{\riskParam}{\FBgain}{\Vmat}}
		[ \InMat^{\MyTRANSPO} \Vmat \DriftMat] 
	}
	\in \mathbb{R}^{  \DimU \DimX   }
	.
\end{align}	
We have ${\ImpTwoF{\SolPair}{\riskParam}}=0$ and ${\Vmat}\succ 0$ if and only if  the corresponding $({\Vmat}\succ 0,{\FBgain})$ is a solution to the WSR  equations \eqref{eq:def_optVmat_inf} and \eqref{eq:def_optFBgain_inf}.
Let  $\DomVmat \subseteq \{ {\Vmat} \in \DomALLVmat | {\Vmat}\succ 0 \} $ and $\DomFBgain \subseteq \DomALLFBgain$ be arbitrarily assigned bounded closed sets.
The corresponding set of  $\SolPair$ is denoted by $\DomSolPair:=\{  \SolPair \in \mathbb{R}^{ (\DimX(\DimX+1)/2) + \DimU \DimX   } |    ({\Vmat},{\FBgain})  \in  \DomVmat \times \DomFBgain   \}$.
We consider solutions  on $\DomVmat \times \DomFBgain$ because this boundedness is reasonable for implementing the controllers.
We introduce the following assumption:

\begin{assumption}\label{ass:weight}
	The PDF $\PDF{\StoParam}$, set $\DomStoParam$,  desired weight $\dWeight$, and set $\DomSolPair$ satisfies the following conditions:		
	\begin{enumerate}
		
		\item[{\AssWeightMSstabilizable}]
		The feedback system \eqref{eq:def_FBsys} is MS stabilizable.
		
		\item[{\AssWeightZeroSolContained}]
		The interior of $\DomSolPair$ contains a solution to the WSR equations \eqref{eq:def_optVmat_inf} and \eqref{eq:def_optFBgain_inf} for $\riskParam=0$.

		\item[{\AssWeightSmooth}]
		There exist an upper bound $\smoothUBriskParam > 0$ and a lower bound $\smoothLBriskParam < 0$ such that
		${\WE{\StoParam}{\riskParam}{\FBgain}{\Vmat}}	[ \StoParam \StoParam^{\MyTRANSPO}]$
		are $C^{2}$ continuous  on an open subset of $\DomALLVmat \times \DomALLFBgain \times  \mathbb{R}$ and this subset contains $ \DomVmat\times\DomFBgain  \times [\smoothLBriskParam,\smoothUBriskParam] $.

	\end{enumerate}

\end{assumption}

\begin{theorem}[Uniqueness and smoothness]\label{thm:uniqueness}
	Suppose that Assumption \ref{ass:weight} holds.
	There exist an upper bound $\uniqueUBriskParam > 0$ and a lower bound $\uniqueLBriskParam < 0$ satisfying the following two statements.
	For each $\riskParam \in [\uniqueLBriskParam,\uniqueUBriskParam]$,
	there exists a unique solution $(\IoptVmat,\IoptFBgain) \in \DomVmat\times\DomFBgain  $ to the WSR  equations \eqref{eq:def_optVmat_inf} and \eqref{eq:def_optFBgain_inf} satisfying $\IoptVmat \succ 0$, provided that the set of solutions are restricted to $\DomVmat\times\DomFBgain $.
	The unique solution $(\IoptVmat,\IoptFBgain)$ is $C^{1}$ continuous in $\riskParam$ on $(\uniqueLBriskParam,\uniqueUBriskParam)$.  
\end{theorem}
\begin{\MyProof}	
	The proof is described in Appendix \ref{pf:uniqueness}.
\end{\MyProof}

\begin{remark}[Contribution of Theorem \ref{thm:uniqueness}]
	Theorem \ref{thm:uniqueness} ensures the uniqueness and smoothness of the solution to the WSR  equations.
	These properties are helpful for solving the WSR  equations in Section \ref{sec_solve_problem_1_solve} and analyzing stability in Section \ref{sec_method2}.
\end{remark}

\subsection{How to solve the WSR equations}\label{sec_solve_problem_1_solve}

We propose two approaches for solving the WSR equations.
The first is to iterate the following WSR difference equations.
\begin{definition}[WSR difference equations]
	Given $\dWeight$, $\riskParam$, $\IiniVmat \in  \DomALLVmat$, and $\IiniFBgain\in \DomALLFBgain$, 
	let us define the WSR difference equations for $\MybT\in \{0,1,2,\dots\}$ as
	\begin{align}	
		\begin{bmatrix}
			{\optVmat{\MybT+1}} 
			\\ 
			{\optFBgain{\MybT+1}} 
		\end{bmatrix}
		&:=
		\begin{bmatrix}
			{\WSRDVmat{\optVmat{\MybT}}{\optFBgain{\MybT}}{\riskParam}}	
			\\
			{\WSRDFBgain{\optVmat{\MybT}}{\optFBgain{\MybT}}{\riskParam}}	
		\end{bmatrix}
		. \label{eq:def_WSRD}		
	\end{align}	
\end{definition}

\begin{theorem}[Solution to the WSR equations]	\label{thm:solution_WSR}
	Given $\dWeight$, $\riskParam$, $\IiniVmat \succeq 0 \in  \DomALLVmat $, and $\IiniFBgain \in \DomALLFBgain$, 
	suppose that 
	\eqref{eq:def_WSRD} is well defined  for every $\MybT\in \{0,1,2,\dots\}$.
	If there exists a pair $(\IestVmat,\IestFBgain)$ that satisfies
	\begin{align}	
		\begin{bmatrix}
			\IestVmat
			\\ 
			\IestFBgain
		\end{bmatrix}		
		=
		\lim_{\MybT \rightarrow \infty} 
		\begin{bmatrix}
			{\optVmat{\MybT}} 
			\\ 
			{\optFBgain{\MybT}} 
		\end{bmatrix}
		,\label{eq:lim_WSRDeq}
	\end{align}
	then $(\IestVmat,\IestFBgain)$ is a solution to the WSR  equations \eqref{eq:def_optVmat_inf} and \eqref{eq:def_optFBgain_inf}  satisfying $\IestVmat\succ 0$,
	provided that ${\WSRDVmat{\wildcard}{\wildcard}{\riskParam}}$ and ${\WSRDFBgain{\wildcard}{\wildcard}{\riskParam}}$ are continuous at $(\IestVmat,\IestFBgain)$. 
\end{theorem}
\begin{\MyProof}
	The proof is described in Appendix \ref{pf:solution_WSR}.
\end{\MyProof}

\begin{remark}[Contribution of Theorem \ref{thm:solution_WSR}]\label{rem:contribution_solution_WSR}
	We obtain a solution to the WSR equations by iterating the WSR difference equations \eqref{eq:def_WSRD} if they converge successfully.
	If Assumption \ref{ass:weight} {\AssWeightMSstabilizable} and $\riskParam=0$ hold, we guarantee that there exists a pair $(\IestVmat,\IestFBgain)$ satisfying \eqref{eq:lim_WSRDeq},
	as described in Lemma \ref{thm:existing_results} {\KoningPiConverge} in Appendix \ref{sec_supporting}.
	The pair $(\IiniVmat,\IiniFBgain)$ is an initial estimate of $(\IestVmat,\IestFBgain)$, which is typically set to $\IiniVmat=0$ and $\IiniFBgain=0$.
	Section \ref{sec:results_existence} demonstrates that the WSR difference equations converge successfully.
\end{remark}

In the second approach, we show that Newton's method \citep[Chapter 5]{Kelley95} can be successfully employed to solve the WSR equations  \eqref{eq:def_optVmat_inf} and \eqref{eq:def_optFBgain_inf} under Assumption \ref{ass:weight}.
For each $\riskParam$, a solution to the WSR equations and its vectorization are explicitly denoted by $(\IoptVmat(\riskParam) \succ 0,\IoptFBgain(\riskParam) )$ and  ${\optSolPair{\riskParam}}$, respectively.
To calculate ${\optSolPair{\riskParam}}$, we apply Newton's method to ${\ImpTwoF{\SolPair}{\riskParam}}$ as follows:
\begin{equation} \label{eq:Newton_iteration}
	\begin{aligned}	
		&
		{\iteSolPair{\IDiteA+1}}={\iteSolPair{\IDiteA}}-
		\Big(
		\frac{\partial {\NonArgImpTwoF} }{\partial \SolPair^{\MyTRANSPO}}
		({\iteSolPair{\IDiteA}},{\riskParam})
		\Big)^{-1}
		{\ImpTwoF{{\iteSolPair{\IDiteA}}}{\riskParam}}
		,
	\end{aligned}
\end{equation}
where the subscript $\IDiteA$ denotes an iteration index.
To analyze convergence of Newton's method, we introduce the following definitions and lemma that are modified versions of {\citep[Assumption 4.3.1, Definition 4.1.1, Theorem 5.1.2]{Kelley95}}.

\begin{definition}[Standard assumptions] \label{def:standard_ass}
	Given an open set $\DomSaSolPair \subset \DomSolPair$ and $\riskParam$, 
	the following conditions are called the {standard assumptions on $(\DomSaSolPair,\riskParam)$}.
	\begin{enumerate}
		\item[{\StandardAssA}] 
		There exists a solution  ${\optSolPair{\riskParam}} \in \DomSaSolPair$ to ${\ImpTwoF{\SolPair}{\riskParam}} = 0$.
		\item[{\StandardAssB}]
		$\partial {\ImpTwoF{\SolPair}{\riskParam}}/ \partial \SolPair^{\MyTRANSPO}
		$ is Lipschitz continuous in $\SolPair$ on $\DomSaSolPair$.
		\item[{\StandardAssC}] 
		$\partial {\ImpTwoF{\optSolPair{\riskParam}}{\riskParam}}/ \partial \SolPair^{\MyTRANSPO}$ is nonsingular.
	\end{enumerate}
\end{definition}

\begin{definition}[q-quadratic property]
	The convergence ${\iteSolPair{\IDiteA}} \to {\optSolPair{\riskParam}}$ is said to be {q-quadratically} if
	${\iteSolPair{\IDiteA}} \to {\optSolPair{\riskParam}}$ and there exists $\NewtonConst>0$ such that
	$\| {\iteSolPair{\IDiteA+1}} - {\optSolPair{\riskParam}} \|
	\leq
	\NewtonConst
	\| {\iteSolPair{\IDiteA}} - {\optSolPair{\riskParam}} \|^{2}$.
\end{definition}

\begin{lemma}[Newton's method] \label{def:original_Newton}
	Given $\DomSaSolPair \subset \DomSolPair$, $\riskParam$, and ${\iteSolPair{0}} \in \DomSolPair$, there exists $\DefBallRadius>0$ such that 	
	if the following conditions (i) and (ii) hold,
	${\iteSolPair{\IDiteA}}$ in \eqref{eq:Newton_iteration} converges {q-quadratically} to a solution ${\optSolPair{\riskParam}}$.
	\begin{enumerate}
		\item 
		The standard assumptions on $(\DomSaSolPair,\riskParam)$ hold.
		\item
		We have
		${\iteSolPair{0}} \in {\BallSolPair{\DefBallRadius}{\riskParam}} \subset \DomSaSolPair$, where
		${\BallSolPair{\DefBallRadius}{\riskParam}}
		:=\{ \SolPair | \|\SolPair-{\optSolPair{\riskParam}}\| < \DefBallRadius \} $ and ${\optSolPair{\riskParam}}$ is unique.
	\end{enumerate}
\end{lemma}

However, guaranteeing the conditions (i) and (ii) in Lemma \ref{def:original_Newton} is not straightforward.
We must ensure that the standard assumptions on $( \DomSaSolPair , \riskParam )$ hold.
Moreover, an initial estimate ${\iteSolPair{0}}$ included in ${\BallSolPair{\DefBallRadius}{\riskParam}}$ must be appropriately determined.
This study overcomes these difficulties and employs Newton's method as follows:

\begin{theorem}[Successful Newton's method] \label{thm:MyNewton}
	Suppose that Assumption \ref{ass:weight} holds.
	Let ${\iteSolPair{0}}:={\optSolPair{0}}$.
	There exist an upper bound $\NewtonUBriskParam>0$ and a lower bound $\NewtonLBriskParam < 0$ such that for every $\riskParam \in (\NewtonLBriskParam,\NewtonUBriskParam)$,
	${\iteSolPair{\IDiteA}}$ in \eqref{eq:Newton_iteration} converges {q-quadratically} to a solution ${\optSolPair{\riskParam}}$.	
\end{theorem}
\begin{\MyProof}
	The proof is described in Appendix \ref{pf:MyNewton}.		
\end{\MyProof}

\begin{remark}[Contribution of Theorem \ref{thm:MyNewton}]
	By virtue of Theorem \ref{thm:MyNewton}, Newton's method is successfully applied to solve the WSR equations \eqref{eq:def_optVmat_inf} and \eqref{eq:def_optFBgain_inf}.
	Theorem \ref{thm:MyNewton} shows that the solution ${\optSolPair{0}}$ with $\riskParam=0$ is a suitable initial estimate.
	It is easy to obtain ${\optSolPair{0}}$ as described in Remarks \ref{rem:special_cases} and \ref{rem:contribution_solution_WSR}.
\end{remark}

\section{Proposed method: solution to Problem 2}\label{sec_method2}

We solve Problem 2: analyzing the WMS stability of the feedback system \eqref{eq:def_FBsys} with applying the controllers derived in the previous section.
Firstly, we discuss the MS stability, which is a special case of the WMS stability.

\begin{proposition}[MS stability discrimination] \label{thm:ms-stablity}
	Given a feedback gain $\FBgain\in \mathbb{R}^{\DimU \times \DimX}$, suppose that $\funcInput(\State)=-\FBgain \State$ holds.
	The feedback system \eqref{eq:def_FBsys} is MS stable if and only if 
	the spectral radius of ${\EliMat{\DimX}}	
	{\Expect{\StoParam}}[{\DescriptionFBMat{\FBgain}} \otimes {\DescriptionFBMat{\FBgain}} ]
	{\DupMat{\DimX}}$, that is, the maximum absolute value of the eigenvalues, is less than 1, where ${\DupMat{\DimX}}$ and ${\EliMat{\DimX}}$ are defined in Section \ref{sec_intro}.
	In addition, the  feedback system is MS stabilizable if and only if there exists $\FBgain$ such that the above spectral radius is less than 1.
\end{proposition}
\begin{\MyProof}	
	The statement follows from the results in \citep[Theorems 1 and 3]{ItoIFACWC20} by removing time-invariant stochastic parameters  because the asymptotic stability of discrete-time linear systems are characterized by the spectral radius being less than one.
\end{\MyProof}

\begin{theorem}[MS stability] \label{thm:ms_stable_existence}
	Suppose that Assumption \ref{ass:weight} holds.
	There exist an upper bound $\StableUBriskParam > 0$ and a lower bound $\StableLBriskParam < 0$ such that
	for all $\riskParam \in ( \StableLBriskParam , \StableUBriskParam)$, 
	the feedback system \eqref{eq:def_FBsys} with applying the optimal controller $\optFBcontroller(\State; {\IoptFBgain(\riskParam)})$ in \eqref{eq:sol_optFBcontroller} is MS stable.
\end{theorem}
\begin{\MyProof}
	The proof is described in Appendix \ref{pf:ms_stable_existence}.	
\end{\MyProof}
\begin{remark}[Contribution of Theorem \ref{thm:ms_stable_existence}]
	Theorem \ref{thm:ms_stable_existence} guarantees that we can design controllers  $\optFBcontroller(\State; {\IoptFBgain(\riskParam)})$ making the feedback system MS stable for each $\riskParam \in ( \StableLBriskParam , \StableUBriskParam)$. 
	By using Proposition \ref{thm:ms-stablity}, we can easily evaluate whether
	the feedback system \eqref{eq:def_FBsys} with $\optFBcontroller(\State; {\IoptFBgain(\riskParam)})$ is MS stable.
	Even if the MS stability is not ensured, 
	decreasing the absolute value of $\riskParam$ enables the design of ${\IoptFBgain(\riskParam)}$ that guarantees the MS stability. 
\end{remark}

Next, we discuss the WMS stability of the feedback system.
Recall that the proposed controllers focus on minimizing the weighted cost function \eqref{eq:def_IHWcostJ}  associated with the weight ${\Weight{}}$.
We show that such weighted control problems are compatible with the WMS stability.

\begin{theorem}[WMS stability]	\label{thm:wms_stable}
	For any $\riskParam$ and any solution $(\IoptVmat(\riskParam),\IoptFBgain(\riskParam)) $ to the WSR equations \eqref{eq:def_optVmat_inf} and \eqref{eq:def_optFBgain_inf} satisfying $\IoptVmat(\riskParam) \succ 0$,
	the feedback system \eqref{eq:def_FBsys} with applying the optimal controller $\optFBcontroller(\State; {\IoptFBgain(\riskParam)})$ in \eqref{eq:sol_optFBcontroller}
	is WMS stable with the following weight functions:
	\begin{equation} \label{eq:def_weight_t}
		\begin{aligned}	
			&
			\forall \MyT ,\;\;
			{\WmsWeight{\MyT}}
			=
			\prod_{\MybT=0}^{\MyT}
			\dWeight(\StoParam_{\MybT}; \riskParam, \optFBcontroller ,  \mincostJ )
			.
		\end{aligned}
	\end{equation}	
\end{theorem}
\begin{\MyProof}
	The proof is described in Appendix \ref{pf:wms_stable}.
\end{\MyProof}

Theorem \ref{thm:wms_stable} indicates that the state $\State_{\MyT}$ converges to zero in terms of WMS.
An interpretation of the WMS stability associated with the MS stability is presented below.

\begin{corollary}[Equivalence of the WMS stability]	\label{thm:characteristic_w}
The WMS stability with ${\WmsWeight{\MyT}}$ in \eqref{eq:def_weight_t} is equivalent to the MS stability with the biased PDF ${\biasPDF{\StoParam}}=\dWeight(\StoParam; \riskParam, \optFBcontroller ,  \mincostJ )\PDF{\StoParam}$, that is,
	\begin{align}		
		\lim_{\MyT\to 0}
		{\Expect{\StoParam_{0:\MyT}}}[ {\WmsWeight{\MyT}} \|\State_{\MyT}\|^{2} ]
		=0
		\Leftrightarrow
		\lim_{\MyT\to 0}
		{\biasE{\StoParam_{0:\MyT}}}[  \|\State_{\MyT}\|^{2} ]
		=0
		,\label{eq:p_to_q1}
	\end{align}
	where ${\biasE{\StoParam}}$ is the expectation with respect to $\StoParam$ that obeys the biased PDF ${\biasPDF{\StoParam}}$ instead of ${\PDF{\StoParam}}$.
\end{corollary}
\begin{\MyProof}
We have 
${\Expect{\StoParam_{0:\MyT}}}[ {\WmsWeight{\MyT}} \|\State_{\MyT}\|^{2} ]
=
{\biasE{\StoParam_{0:\MyT}}}[  \|\State_{\MyT}\|^{2} ]
$ in the proof of Theorem \ref{thm:wms_stable}, 
 implying \eqref{eq:p_to_q1}.
\end{\MyProof}
\begin{remark}[Contribution of Corollary \ref{thm:characteristic_w}]
	Owing to \eqref{eq:p_to_q1}, the WMS stability guarantees the MS stability with the PDF  ${\biasPDF{\StoParam}}$ biased by the desired weight $\dWeight(\StoParam, \riskParam, \optFBcontroller ,  \mincostJ )$.
	This bias handles the importance of each value of $\StoParam_{\MyT}$.
	For example, the RSL and RRSL controllers for $\riskParam>0$ add a bias that mitigates worse cases among various control results.
\end{remark}

\section{Numerical example} \label{sec_numerical_example}

This section presents numerical examples to demonstrate the effectiveness of the proposed controllers in terms of 
convergence, robustness, stability, and control performance.
We compare the RRSL controllers with a risk neutral (RN) controller ($\riskParam \to 0$) and RSL controllers.

\subsection{Plant system and setting} \label{sec:plant}

Let us consider the linear system \eqref{eq:def_sys}, where $\DriftMat_{\MyT}$ and $\InMat_{\MyT}$ obey the normal and Laplace distributions, respectively.
The mean and covariance of the matrices are set as follows:
\begin{align}
	\Expect{\StoParam}[\DriftMat] 
	&= \begin{bmatrix}
		0.97&-0.03\\
		0.1&1.03
	\end{bmatrix}	
	,\;\;
	\Expect{\StoParam}[\InMat] 
	= \begin{bmatrix}
		0.005\\0.01
	\end{bmatrix}	,\label{eq:def_AB}
	\\
	\Covariance{\StoParam}[\StoParam]	
	& 
	=({\DIAG{  {\Expect{\StoParam}}[\StoParam] }} /10 )^{2}
	, \label{eq:def_ABsigma}
\end{align}
where 
$\DIAG{   {\Expect{\StoParam}}[\StoParam] }$ denotes the diagonal matrix such that its diagonal components are $ {\Expect{\StoParam}}[\StoParam]$. 
The parameters of the cost in \eqref{eq:def_IHWcostJ} are set as $\xCostMat = 3 {\Identity{2}}$ and $\uCostMat =1$.

\subsection{Existence of the proposed controllers}  \label{sec:results_existence} 

The convergence of the WSR difference equations \eqref{eq:def_WSRD} is evaluated.
We adopt the RRSL control with the weight in Example \ref{rem:RRSL}.
The PDF of $\RefX$ in \eqref{eq:def_predcostJ} is set such that ${\Expect{\RefX}}[ \RefX \RefX^{\MyTRANSPO}] = {\Identity{2}}$ holds.
The other parameters are set as $\paramWeightA=10$ and $\paramWeightB=11$.
The expectations $\Expect{\StoParam}[\dots]$ are calculated using the MC approximation with 10,000 random samples. 
Figure \ref{fig:Pmat_conv} shows the sequences of ${\optVmat{\MybT}}$ and ${\optFBgain{\MybT}}$ with  the initial values $\IiniVmat=0$ and $\IiniFBgain=0$. 
The sequences converged sufficiently.
This convergence implies that a solution to Problem 1 and the corresponding optimal controller were numerically obtained.
In the following, we regard  ${\optVmat{\MybT}}$ and ${\optFBgain{\MybT}}$ for $\MybT=300$ as ${\IoptVmat}$ and ${\IoptFBgain}$, respectively.

\definecolor{myPlotR}{rgb}{ 1,0,0 }
\definecolor{myPlotG}{rgb}{ 0,0.5,0 }
\definecolor{myPlotB}{rgb}{ 0,0,1 }

\begin{figure}[t]\MyFigRePos{.55\linewidth}{	
	\begin{minipage}[b]{.5\linewidth}
		\centering			
		{\scriptsize		
			\hfill
			\begin{tabular}{|l|}
				\hline
				\textcolor{myPlotR}{\textbf{---}}:${\El{\optVmat{\MybT}}{1,1}}$
				,
				\textcolor{myPlotG}{\textbf{$- \cdot -$}}:${\El{\optVmat{\MybT}}{2,1}}$
				\\
				\textcolor{myPlotB}{\textbf{$-\; -$}}:${\El{\optVmat{\MybT}}{2,2}}$
				\\\hline				
			\end{tabular}
		}%
		\begin{tikzpicture}[scale=1.0]
		\begin{axis}[axis y line=none, axis x line=none
		,xtick=\empty, ytick=\empty
		,xmin=-10,xmax=10,ymin=-10,ymax=10
		,width=2.25 in,height=1.55 in
		]
		
		\node at (0,0) { 
			\includegraphics[width=0.95\linewidth]{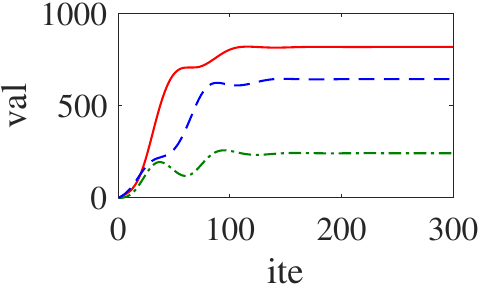}
		};
		
		\filldraw[draw=white, fill=white] (-15,-15) -- (-15,15) -- (-8,15) -- (-8,-15);	
		\filldraw[draw=white, fill=white] (-15,-15) -- (-15,-7.6) -- (15,-7.6) -- (15,-15);	
		\node at (-9.2,2.5) {\rotatebox{90}{\small	
				{\shortstack{${\El{\optVmat{\MybT}}{\IDEl,\IDbEl}}$}} 
		}};
		\node at (1,-8.7) {\rotatebox{0}{\small	
				{\shortstack{The number $\MybT$ of iterations}} 
		}};
		
		\end{axis}			
		\end{tikzpicture}		
		\subcaption{Transition of ${\optVmat{\MybT}}$}\label{fig:Pmat_conv_RR}
	\end{minipage}%
	\begin{minipage}[b]{.5\linewidth}
		{\scriptsize
			\hfill	
			\begin{tabular}{|l|}
				\hline
				\textcolor{myPlotR}{\textbf{---}}:${\El{\optFBgain{\MybT}}{1,1}}$
				\\
				\textcolor{myPlotB}{\textbf{$-\; -$}}:${\El{\optFBgain{\MybT}}{1,2}}$
				\\\hline				
			\end{tabular}
		}%
		\centering
		\begin{tikzpicture}[scale=1.0]
		\begin{axis}[axis y line=none, axis x line=none
		,xtick=\empty, ytick=\empty
		,xmin=-10,xmax=10,ymin=-10,ymax=10
		,width=2.25 in,height=1.55 in
		]
		
		\node at (0,0) { 
			\includegraphics[width=0.95\linewidth]{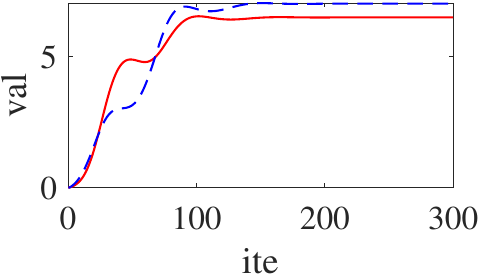}
		};
		
		\filldraw[draw=white, fill=white] (-15,-15) -- (-15,15) -- (-8.3,15) -- (-8.3,-15);	
		\filldraw[draw=white, fill=white] (-15,-15) -- (-15,-7.6) -- (15,-7.6) -- (15,-15);	
		\node at (-9.2,2.5) {\rotatebox{90}{\small	
				{\shortstack{${\El{\optFBgain{\MybT}}{\IDEl,\IDbEl}}$}} 
		}};
		\node at (0,-8.7) {\rotatebox{0}{\small	
				{\shortstack{The number $\MybT$ of iterations}} 
		}};
		
		\end{axis}			
		\end{tikzpicture}		
		\subcaption{Transition of ${\optFBgain{\MybT}}$}\label{fig:FBgain_conv_RR}
	\end{minipage}%
	\caption{Convergence of $({\optVmat{\MybT}} ,	{\optFBgain{\MybT}} )$ in the WSR difference equations for $\riskParam=1$.}\label{fig:Pmat_conv}
}\end{figure}

\subsection{Robustness of the RRSL controllers}\label{sec:robustness}

We evaluate the robustness of the RRSL controllers with respect to random samples when the MC method approximates the expectations $\Expect{\StoParam}[\dots]$.
The MC method is promising for approximating the expectations if they are not obtained in an analytical manner. 
The robustness is compared between the RRSL and RSL controllers.
For several sensitivity parameters $\riskParam$, both the feedback gains ${\optFBgain{\MybT}}$ for $\MybT=300$ are designed 100 times by changing random seeds.
The number of random samples for the MC method are set to 10,000 for each design.
The means and standard deviations of the designed gains are presented by markers and error bars, respectively, in Fig. \ref{fig:variations}.
The standard deviation of the gain for each RRSL controller was less than that of the RSL controller close to the RRSL controller.
This result shows that  RRSL control is superior to  RSL control in terms of the robustness to  random samples.
The means of the designed gains are used for evaluations in the following Sections \ref{sec:results_stability} and \ref{sec:results_performance}.

\definecolor{myColorPlotA}{rgb}{ 0,0,0 }
\definecolor{myColorPlotB}{rgb}{ 0,0,1 }
\definecolor{myColorPlotC}{rgb}{ 1,0,0 }
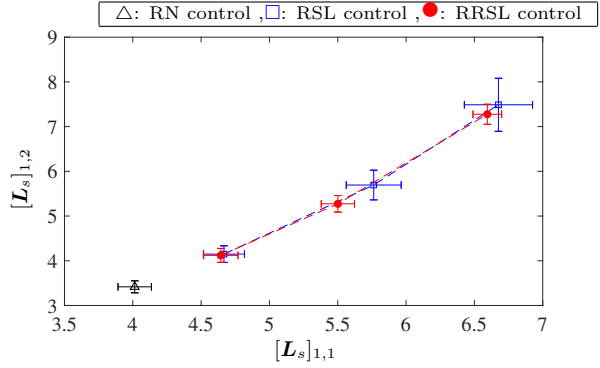
\begin{figure}	\MyFigRePos{.54\linewidth}{		
	{\scriptsize
		\hfill
		\begin{tabular}{| l |}
			\hline
			$\textcolor{myColorPlotA}{\triangle}$: RN control
			,
			$\textcolor{myColorPlotB}{\square}$: RSL control
			,
			{\Large$\textcolor{myColorPlotC}{ \bullet}$}: RRSL control 
			\\
			\hline		
		\end{tabular}	
	}
	\centering	
	\begin{tikzpicture}[scale=0.9]
	\begin{axis}[axis y line=none, axis x line=none
	,xtick=\empty, ytick=\empty
	,xmin=-10,xmax=10,ymin=-10,ymax=10
	,width=3.9 in,height=2.6 in
	]
	
	\node at (0,0) { 
		\includegraphics[width=0.94\linewidth]{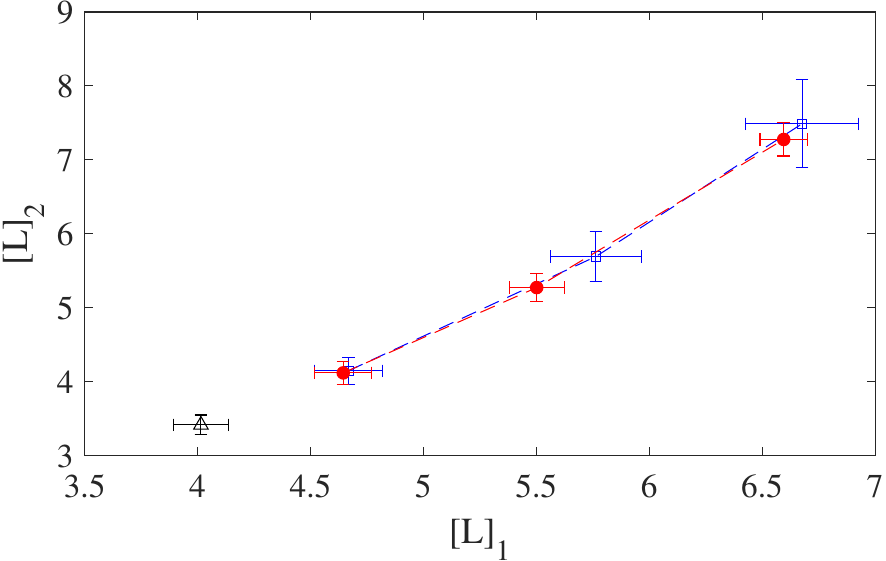}
	};
	
	\filldraw[draw=white, fill=white] (-15,-15) -- (-15,15) -- (-8.4,15) -- (-8.4,-15);	
	\filldraw[draw=white, fill=white] (-15,-15) -- (-15,-8.0) -- (15,-8.0) -- (15,-15);	
	
	\node at (-9.2,1.4) {\rotatebox{90}{\small	
			{\shortstack{ ${\El{\optFBgain{\MybT}}{1,2}}$ }} 
	}};
	\node at (1,-8.9) {\rotatebox{0}{\small	
			{\shortstack{ ${\El{\optFBgain{\MybT}}{1,1}}$ }} 
	}};

	\node at (-4.0,6.5) {
		
	};
	
	\end{axis}			
	\end{tikzpicture}	
	
	\caption{Variations in three types of controllers for several sensitivity parameters $\riskParam$  with $\MybT=300$.
		For the RSL control, $\riskParam$ was set to $0.0005$ (left), $0.001$ (center), and $0.00125$ (right).
		For the RRSL control, $\riskParam$ was set to $0.2$ (left), $0.5$ (center), and $1.0$ (right).
	}	
	\label{fig:variations}
}\end{figure}

\subsection{MS stability}\label{sec:results_stability}

\begin{figure}[t!]\MyFigRePos{.55\linewidth}{			
	\centering	
	\begin{tikzpicture}[scale=0.9]
	\begin{axis}[axis y line=none, axis x line=none
	,xtick=\empty, ytick=\empty
	,xmin=-10,xmax=10,ymin=-10,ymax=10
	,width=2.6 in,height=1.65 in
	]

	\node at (0,0) { 
		\includegraphics[width=0.5\linewidth]{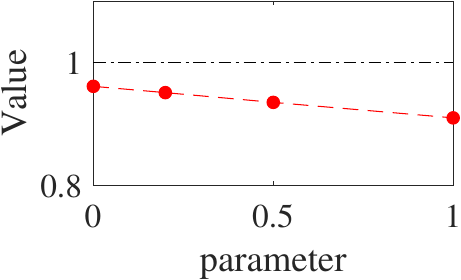}
	};

	\filldraw[draw=white, fill=white] (-15,-15) -- (-15,-7.6) -- (15,-7.6) -- (15,-15);	
	\node at (-8.4,1.4) {\rotatebox{90}{\small	
			\colorbox{white}{\shortstack{Spectral radius}} 
	}};
	\node at (1,-8.7) {\rotatebox{0}{\small	
			{\shortstack{Sensitivity parameter $\riskParam$ }} 
	}};
	
	\end{axis}			
	\end{tikzpicture}	
	\caption{Spectral radius of $
			{\EliMat{\DimX}}	
			{\Expect{\StoParam}}[{\DescriptionFBMat{\optFBgain{\MybT}}} \otimes {\DescriptionFBMat{\optFBgain{\MybT}}} ]
			{\DupMat{\DimX}}
		$ for $\MybT=300$.}	
	\label{fig:RS_PsikronPsi}
}\end{figure}

We evaluate the MS stability of the feedback systems with applying the RRSL controllers.
As described in Proposition \ref{thm:ms-stablity}, the MS stability is discriminated by the spectral radius.
Figure \ref{fig:RS_PsikronPsi} shows the spectral radii for different values of the sensitivity parameter $\riskParam$.
Because the spectral radius was less than 1 for every selected $\riskParam$, Proposition \ref{thm:ms-stablity} theoretically guarantees the MS stability of the feedback system.

\begin{figure}[t]\MyFigRePos{.55\linewidth}{
	\begin{minipage}[b]{.5\linewidth}
		\centering
		\includegraphics[width=0.95\linewidth]{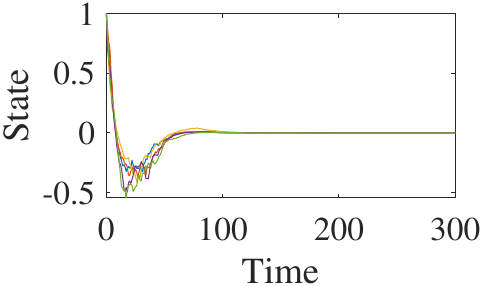}
		\subcaption{First component $\El{\State_{\MyT}}{1}$}\label{fig:RRS_x1}
	\end{minipage}%
	\begin{minipage}[b]{.5\linewidth}
		\centering
		\includegraphics[width=0.95\linewidth]{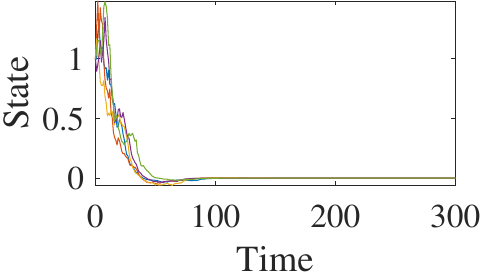}
		\subcaption{Second component $\El{\State_{\MyT}}{2}$}\label{fig:RRS_x2}
	\end{minipage}%
	\caption{State transitions for the proposed RRSL control ($\riskParam=1$).}\label{fig:RSe_stability}
}\end{figure}

Figure \ref{fig:RSe_stability} shows examples of the state transitions of the feedback system \eqref{eq:def_FBsys} when applying the RRSL controller, where multiple trajectories represent different trials.
The initial state is set to $\State_{0}:=[1,1]^{\MyTRANSPO}$.
We can see that the state $\State_{\MyT}$ converged to zero. 	
Therefore, the feedback system was numerically shown to be MS stable. 

\subsection{Control performance}\label{sec:results_performance}

The next evaluation focuses on the control performance.
We compared the RRSL controllers to an existing RN controller ($\riskParam=0$) that is equivalent to a standard stochastic optimal controller  \citep{Koning82}.
The expected quadratic cost ${\Expect{\StoParam_{0:\terminalT}}}[
\sum_{\MyT=0}^{\terminalT}   ( \State_{\MyT}^{\MyTRANSPO} \xCostMat \State_{\MyT} + \Input_{\MyT}^{\MyTRANSPO} \uCostMat \Input_{\MyT}  ) 
]$ is used as the control performance, where $\terminalT=300$ and $\State_{0}=[1,1]^{\MyTRANSPO}$.
We simulate the RRSL control $(\riskParam > 0)$ and the baseline RN control $(\riskParam=0)$ in $100,000$ trials for each $\riskParam$.

Figure \ref{fig:percentiles} compares the average values of worse costs over the $100,000$ trials.
For each $\MyPercentWorse$, the average value denotes the mean value of the costs ranked in the worst $\MyPercentWorse \%$ in all trials.
For example, the average value for $\MyPercentWorse=20$ is calculated as the mean value of the costs ranked in the worst $20,000$ trials.
Large values of the sensitivity parameter $\riskParam$ significantly suppressed the poor results.
These results indicate that the proposed RRSL controller successfully handles the risks of the control results.

\definecolor{myColorBarA}{rgb}{ 0,0,0 }
\definecolor{myColorBarB}{rgb}{ 1.0,0.0,0.0 }
\definecolor{myColorBarC}{rgb}{ 1.0,0.4,0.4 }
\definecolor{myColorBarD}{rgb}{ 1.0,0.8,0.8 }
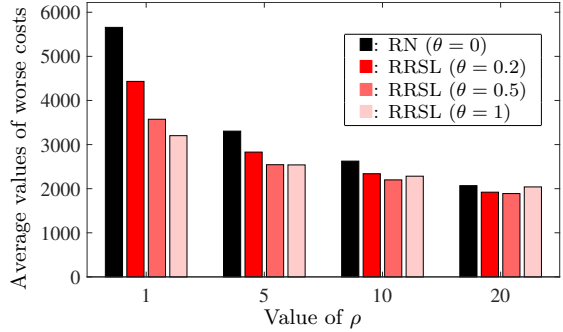
\begin{figure}[t!]\MyFigRePos{.55\linewidth}{			
	\centering	
	\begin{tikzpicture}[scale=0.9]
	\begin{axis}[axis y line=none, axis x line=none
	,xtick=\empty, ytick=\empty
	,xmin=-10,xmax=10,ymin=-10,ymax=10
	,width=3.9 in,height=2.55 in
	]
	
	\node at (0,0) { 
		\includegraphics[width=0.94\linewidth]{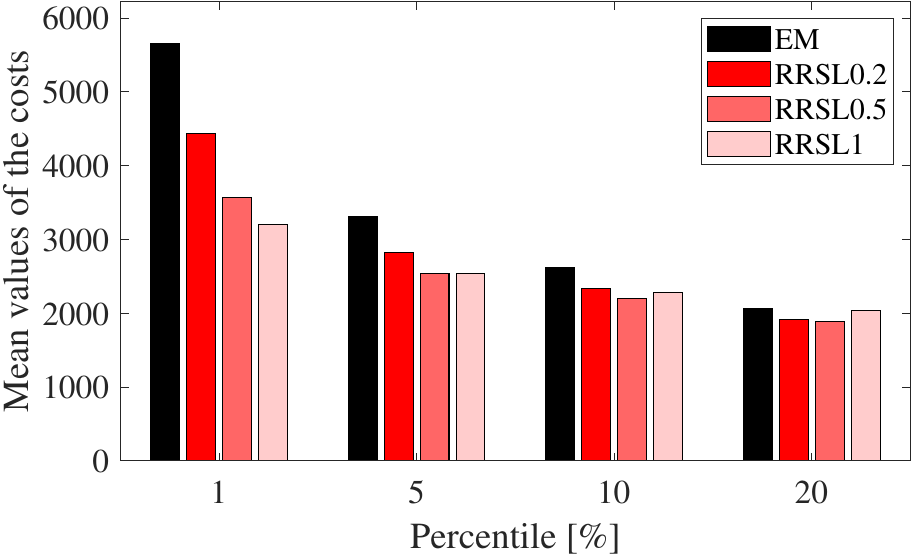}
	};

	\filldraw[draw=white, fill=white] (-15,-15) -- (-15,15) -- (-8.8,15) -- (-8.8,-15);	
	\filldraw[draw=white, fill=white] (-15,-15) -- (-15,-8.2) -- (15,-8.2) -- (15,-15);

	\node at (1,-9.1) {\rotatebox{0}{\small	
			{\shortstack{$\qquad$ Value of $\MyPercentWorse$ $\qquad$ }} 
	}};
	\node at (-9.5,1.4) {\rotatebox{90}{\small	
			{\shortstack{Average values of worse costs}} 
	}};

	\fill[fill=white] (4,2) rectangle (9.4,9.45);
	\node at (5.5,5.3) {
		{\footnotesize
			\begin{tabular}{| l |}
			\hline
			$\textcolor{myColorBarA}{\blacksquare}$: RN ($\riskParam=0$)
			\\
			$\textcolor{myColorBarB}{\blacksquare}$: RRSL ($\riskParam=0.2$)
			\\
			$\textcolor{myColorBarC}{\blacksquare}$: RRSL ($\riskParam=0.5$)   
			\\
			$\textcolor{myColorBarD}{\blacksquare}$: RRSL ($\riskParam=1$)
			\\
			\hline			
			\end{tabular}	
		}
	};
	
	\end{axis}			
	\end{tikzpicture}	
	\caption{Average values of worse performance results.
		For each $\MyPercentWorse$, each average cost denotes the mean value of costs ranked in the worst $\MyPercentWorse \%$ of all the trials.}\label{fig:percentiles}
}\end{figure}

\section{Conclusion}\label{sec_conclusion}

This paper presented a general framework for designing multiple optimal controllers for linear systems with i.i.d. stochastic parameters.
They include stochastic optimal, RSL, and the proposed RRSL controllers.
The proposed WSR equations are powerful tools for deriving the controllers.
The WSR equations are solved via two approaches: using the WSR difference equations and Newton's method.
The stability of the feedback systems was analyzed in terms of MS and WMS.

The proposed general theory has the potential to provide even more types of optimal controllers, by designing weight functions.
A further challenge is to find novel optimal controllers using the WSR equations.
Additionally, extending the proposed method to other types of stochastic parameters broadens its applicability.


\appendix
\section*{Appendix}
\addcontentsline{toc}{section}{Appendix}

\section{Supporting results}\label{sec_supporting}

Let ${\biasE{\StoParam}}$ denote the expectation with respect to $\StoParam$ that follows another PDF ${\biasPDF{\StoParam}}$ instead of ${\PDF{\StoParam}}$.
For any $\FBgain$, we define
the linear operator ${\Ltrans{}{\FBgain}{\NotationSymMat}{}} 
:= {\biasE{\StoParam}}[ {\DescriptionFBMat{\FBgain}}{}^{\MyTRANSPO} \NotationSymMat {\DescriptionFBMat{\FBgain}}]$,
matrix ${\kronLtrans{\FBgain}}:={\biasE{\StoParam}}[ $
${\DescriptionFBMat{\FBgain}}{}^{\MyTRANSPO} \otimes {\DescriptionFBMat{\FBgain}}{}^{\MyTRANSPO}]$, and
compression operator ${\Elimi{\NotationSquMat}}:={\EliMat{\DimX}} \NotationSquMat {\DupMat{\DimX}}$ \citep[Definition 3]{ItoIFACWC20}.

\begin{definition}[Stable linear transformations]\label{def:stable_trans}
	Given $\FBgain$, the function ${\NonArgLtrans{\FBgain}{}}$ is said to be stable if 
	the spectral radius of ${\Elimi{ {\kronLtrans{\FBgain}} }}$, that is, the maximum absolute value of the eigenvalues, is less than 1.
\end{definition}

\begin{lemma}[Equivalence of stable properties]\label{thm:def_stable2}
	The following conditions are equivalent.
	\begin{enumerate}
		\item 
		The operator ${\NonArgLtrans{\FBgain}{}}$ is stable.
		\item
		For some $\NotationSymBMat \succ 0 \in {\SetSymMat{\DimX}}$,
		 $\lim_{\MybT \rightarrow \infty } {\Ltrans{}{\FBgain}{\NotationSymBMat}{\MybT}} $ $ =0$.
		\item
		For any $\NotationSymMat \succeq 0 \in {\SetSymMat{\DimX}}$,
	 $\lim_{\MybT \rightarrow \infty } {\Ltrans{}{\FBgain}{\NotationSymMat}{\MybT}} =0$.
		
	\end{enumerate}
	
\end{lemma}
\begin{\MyProof}	
	Firstly, we prove (i) $\Rightarrow$ (ii) similarly to \citep[Appendix A]{ItoIFACWC20}.
	We have
	$\VECH{ {\Ltrans{}{\FBgain}{\NotationSymMat}{}} } 
	= 
	{\EliMat{\DimX}} \VEC{ {\Ltrans{}{\FBgain}{\NotationSymMat}{}} } 
	=
	{\EliMat{\DimX}} {\kronLtrans{\FBgain}} \VEC{ \NotationSymMat } 
	=
	{\Elimi{\kronLtrans{\FBgain}}} \VECH{\NotationSymMat}$
	and thus
	$\VECH{ {\Ltrans{}{\FBgain}{\NotationSymMat}{\MybT}} }=  {\Elimi{\kronLtrans{\FBgain}}} $ $ \VECH{ {\Ltrans{}{\FBgain}{\NotationSymMat}{\MybT-1}} }
	= \dots = {\Elimi{\kronLtrans{\FBgain}}}^{\MybT} {\VECH{\NotationSymMat}} $
	for any symmetric $\NotationSymMat$.
	Thus, (i) implies (ii) because 
	$
	\VECH{\Ltrans{}{\FBgain}{\NotationSymBMat}{\MybT}}
	=
	{\Elimi{\kronLtrans{\FBgain}}}^{\MybT}  \VECH{\NotationSymBMat} \to 0$ holds as $\MybT \to \infty$.

	Next, we prove (ii) $\Rightarrow$ (iii) in a manner similar to \citep[the proof of Lemma 2.1]{Koning82}.
	For any $\NotationSymMat \succeq 0$ and $\NotationSymBMat \succ 0$, there exists $\COEFpsdMatCompare \in \mathbb{R}$ such that $0 \preceq \NotationSymMat \preceq \COEFpsdMatCompare \NotationSymBMat$.
	Because of the monotonicity of ${\NonArgLtrans{\FBgain}{\MybT}}$, 
	we obtain 
	$0 = 
	{\Ltrans{}{\FBgain}{0 }{\MybT}} 
	\preceq 
	{\Ltrans{}{\FBgain}{\NotationSymMat}{\MybT}} 
	\preceq
{\Ltrans{}{\FBgain}{\COEFpsdMatCompare \NotationSymBMat}{\MybT}}  
	= \COEFpsdMatCompare
	 {\Ltrans{}{\FBgain}{\NotationSymBMat}{\MybT}} 
	\to 0$ as  $\MybT \to \infty$. 
	This implies that (ii) $\Rightarrow$ (iii).

	Next, we prove (iii) $\Rightarrow$ (i) similar to \citep[Appendix C]{ItoIFACWC20}.
	Define $\NotationSymMat_{\IDEl,\IDbEl} :=( {\unitVec{\IDEl}} + {\unitVec{\IDbEl}} ) ( {\unitVec{\IDEl}} + {\unitVec{\IDbEl}} )^{\MyTRANSPO} \succeq 0$ for $\IDEl \in \{1,\dots, \DimX\}$ and $\IDbEl \in \{\IDEl,\dots, \DimX\}$,
	where ${\unitVec{\IDNotation}} \in \mathbb{R}^{\DimX}$ satisfies ${\El{\unitVec{\IDNotation}}{\IDNotation}}=1$ and the other components are zero.
	From (iii), for any $\VECH{\NotationSymMat_{\IDEl,\IDbEl}} \in \mathbb{R}^{  \DimX(\DimX +1)/2}$, 
	we have
	$
	{\Elimi{\kronLtrans{\FBgain}}}^{\MybT} \VECH{\NotationSymMat_{\IDEl,\IDbEl}} 
	=  
	\VECH{{\Ltrans{}{\FBgain}{\NotationSymMat_{\IDEl,\IDbEl}}{\MybT}}} 
	\to 0$ as  $\MybT \to \infty$.
	Because the basis on $\mathbb{R}^{  \DimX(\DimX +1)/2}$ consists of linear combinations of $\VECH{\NotationSymMat_{\IDEl,\IDbEl}}$,
	for any $\NotationVec  \in \mathbb{R}^{  \DimX(\DimX +1)/2}$, 
we have $
	 {\Elimi{\kronLtrans{\FBgain}}}^{\MybT} \NotationVec \to 0$, which leads to
	  $
	  {\Elimi{\kronLtrans{\FBgain}}}^{\MybT} \to 0$. 
	Because the spectral radius of ${\Elimi{\kronLtrans{\FBgain}}}$ is less than 1, (i) holds.
	This completes the proof.
\end{\MyProof}

We review an existing stochastic optimal control problem \citep{Koning82}.
\begin{definition}[Stochastic optimal control problem]\label{def:SSOCP}
	Given a PDF $\biasPDF{\StoParam}$ of $\StoParam$, consider the system \eqref{eq:def_sys} with $\StoParam_{\MyT}$ that obeys ${\biasPDF{\StoParam}}$ instead of ${\PDF{\StoParam}}$. 
	Consider the cost function
	${\IHWcostJ{\Input_{0:\infty}}{\State_{0}}{\riskParam}} \SSOCPsetting$,
	where $ \SSOCPsetting$ denotes the following three conditions: 
	1) $\dWeight={\Weight{}}=1$ holds for all $\StoParam$,
	2) the expectations are taken with a biased PDF ${\biasPDF{\StoParam}}$, that is, ${\Expect{\StoParam}}$ is replaced with ${\biasE{\StoParam}}$, 
	and 
	3) the controller is not restricted as a time-invariant feedback controller, that is, the argument $\funcInput(\wildcard)$ can be replaced with $\Input_{0:\infty}$.
	With these settings, Problem 1 is said to be the standard stochastic optimal control problem (SSOCP) with $\biasPDF{\StoParam}$.
\end{definition}

\begin{lemma}[{\citep{Koning82}}]\label{thm:existing_results}
	If the SSOCP with $\biasPDF{\StoParam}$ is considered, the following results hold.
	\begin{enumerate}

		\item[{\KoningTransition}]
		(Lemma 3.1)
		For any $\NotationSymMat \in  {\SetSymMat{\DimX}}$ and any $\FBgain$, the feedback system \eqref{eq:def_FBsys} with $\funcInput(\State_{\MyT})=-\FBgain\State_{\MyT}$ satisfies	${\biasE{\StoParam_{0:\MyT-1}}}[\State_{\MyT}^{\MyTRANSPO} {\NotationSymMat} \State_{\MyT}]	
		=\State_{0}^{\MyTRANSPO} {\Ltrans{}{\FBgain}{\NotationSymMat}{\MyT}} \State_{0}$.	
		
		\item
		[{\KoningStable}]  
		(Theorem 3.2)
		Given $\FBgain$, the feedback system \eqref{eq:def_FBsys} with $\funcInput(\State_{\MyT})=-\FBgain\State_{\MyT}$ is MS stable 
		if and only if
		${\NonArgLtrans{\FBgain}{}}$ is stable.

		\item[{\KoningStabilizable}]   
		(Theorem 3.3)
		The system \eqref{eq:def_FBsys} is MS stabilizable 
		if and only if
		there exists $\FBgain$ such that 	${\NonArgLtrans{\FBgain}{}}$  is stable.
		
		\item[{\KoningTransStable}]   
		(Lemma 2.2)
	Given $\NotationSymBMat \succ 0 \in {\SetSymMat{\DimX}}$, there exists a solution $\NotationSymAMat \succ 0 \in {\SetSymMat{\DimX}}$ to $\NotationSymAMat=  {\Ltrans{}{\FBgain}{\NotationSymAMat}{}}  + \NotationSymBMat$ if and only if ${\NonArgLtrans{\FBgain}{}}$ is stable.

		\item[{\KoningPiConverge}]  
		(Theorems 4.3 and 5.1) 
		Suppose that the feedback system \eqref{eq:def_FBsys} is MS stabilizable.
		Then, there exists
		$\KoIoptVmat=\lim_{\MybT\to \infty} {\biasWSRDVmat[\MybT]{0}}$, and it is the minimal nonnegative definite solution to 
		$\Vmat = {\biasWSRDVmat{\Vmat}}$:
		\begin{align}	
			{\biasWSRDVmat{\Vmat}}
			&:=
			{\biasE{\StoParam}}
			[ \DriftMat^{\MyTRANSPO} \Vmat \DriftMat] 
			+
			\xCostMat
			-
			{\biasE{\StoParam}}
			[ \DriftMat^{\MyTRANSPO} \Vmat \InMat] 
			{\biasWSRDFBgain{\Vmat}}
			\nonumber\\&
			=
			{\Ltrans{}{{\biasWSRDFBgain{\Vmat}}}{\Vmat}{}}
			+  {\biasWSRDFBgain{\Vmat}}^{\MyTRANSPO} \uCostMat {\biasWSRDFBgain{\Vmat}}   
			+ \xCostMat
			,\label{eq:def_optVmat_inf_bias}
			\\
			{\biasWSRDFBgain{\Vmat}}
			&:=
			{\biasE{\StoParam}}
			[ \InMat^{\MyTRANSPO} \Vmat \InMat  +  \uCostMat]^{-1}
			{\biasE{\StoParam}}
			[ \InMat^{\MyTRANSPO} \Vmat \DriftMat] 
			. \label{eq:def_WSRDFBgain_bias}
		\end{align}	
		\item[{\KoningIHopt}] 
		(Theorem 5.2)
		Suppose that there exists $\KoIoptVmat=\lim_{\MybT\to \infty} {\biasWSRDVmat[\MybT]{0}}$.
	Then, we have
		$\State_{0}^{\MyTRANSPO} \KoIoptVmat \State_{0}
		=
		\min_{\Input_{0:\infty}}$ $ {\IHWcostJ{{\Input_{0:\infty}}}{\State_{0}}{\riskParam}} \SSOCPsetting
		$, and the optimal input is  $\Input_{\MyT}=\KoIoptFBgain\State_{\MyT}$ with  $\KoIoptFBgain:={\biasWSRDFBgain{\KoIoptVmat}}$.

		\item[{\KoningUnique}]  
		(Theorem 5.3)
		Suppose that there exists $\KoIoptVmat=\lim_{\MybT\to \infty} {\biasWSRDVmat[\MybT]{0}}$.
		Then, $\KoIoptVmat$ is positive definite and the unique  solution to $\Vmat = {\biasWSRDVmat{\Vmat}}$. 
		In addition, ${\NonArgLtrans{\KoIoptFBgain}{}}$  is stable.

	\end{enumerate}
\end{lemma}
\begin{remark}
	The definitions of the stable and stabilizable properties of ${\NonArgLtrans{\FBgain}{}}$ in this study are different from those in \citep{Koning82}.
	Regardless of this difference, Lemma \ref{thm:existing_results} is proven in the same manner as the original version,
	using Lemma \ref{thm:def_stable2}.
	For example, 
	similarly to \citep[Lemma 2.2 (a)]{Koning82},
	if ${\NonArgLtrans{\FBgain}{}}$ is stable,
	there exists $\NotationSymAMat= \sum_{\MybT=0}^{\infty} {\Ltrans{}{\FBgain}{\NotationSymBMat}{\MybT}}$
	because 
	$ \VECH{\NotationSymAMat}
	=\sum_{\MybT=0}^{\infty}\VECH{\Ltrans{}{\FBgain}{\NotationSymBMat}{\MybT}}
	=\sum_{\MybT=0}^{\infty}{\Elimi{\kronLtrans{\FBgain}}}^{\MybT}  \VECH{\NotationSymBMat}
	$ and thus
	$(\Identity{} - {\Elimi{\kronLtrans{\FBgain}}}  ) \VECH{\NotationSymAMat} = {\VECH{\NotationSymBMat}} $
	hold with the nonsingular $(\Identity{} - {\Elimi{\kronLtrans{\FBgain}}}  )$.
	In addition, 
	$(\Identity{} - {\Elimi{\kronLtrans{\FBgain}}}  ) \VECH{\NotationSymAMat} = {\VECH{\NotationSymBMat}} $
	is equivalent to $\NotationSymAMat=  {\Ltrans{}{\FBgain}{\NotationSymAMat}{}}  + \NotationSymBMat$.
	Furthermore, our definitions are advantageous because the spectral radius of ${\Elimi{\kronLtrans{\FBgain}}} $ is easy to evaluate. 
Lemma \ref{thm:existing_results} {\KoningStable} and {\KoningStabilizable} are equivalent to Proposition \ref{thm:ms-stablity} if ${\biasPDF{\StoParam}}={\PDF{\StoParam}}$ holds.
\end{remark}

\section{Proof of Theorem \ref{thm:solution_problem_1}} \label{pf:solution_problem_1}

Firstly, we show that some statements in Lemma \ref{thm:existing_results} in Appendix \ref{sec_supporting} can be utilized because there exists a solution $(\IoptVmat,\IoptFBgain)$ to the WSR equations \eqref{eq:def_optVmat_inf} and \eqref{eq:def_optFBgain_inf}. 
Let us consider the system \eqref{eq:def_sys} with $\StoParam_{\MyT}$ that obeys the biased PDF ${\biasPDF{\StoParam}}=
\dWeight(\StoParam; \riskParam, \estfuncInput(\wildcard;\IoptFBgain)  ,  \estfunccostJ(\wildcard;\IoptVmat)  )
\PDF{\StoParam}$ called the biased system.
Because of Assumption \ref{ass:desired_weight}, ${\biasPDF{\StoParam}}$ is a continuous PDF satisfying $\int_{\DomStoParam} {\biasPDF{\StoParam}} \mathrm{d} \StoParam=1$ and ${\biasPDF{\StoParam}}\geq 0$.
Then, the WSR equations with $(\IoptVmat,\IoptFBgain)$ reduce to
$\IoptVmat = {\biasWSRDVmat{\IoptVmat}}$ and  $\IoptFBgain:={\biasWSRDFBgain{\IoptVmat}}$ with 
 \eqref{eq:def_optVmat_inf_bias} and \eqref{eq:def_WSRDFBgain_bias}  because the weighted expectation ${\WE{\StoParam}{\riskParam}{\IoptFBgain}{\IoptVmat}}$ is replaced with the expectation ${\biasE{\StoParam}}$.
Because $\IoptVmat$ is a positive definite solution to \eqref{eq:def_optVmat_inf_bias}, the biased system is MS stabilizable by Lemma \ref{thm:existing_results} {\KoningTransStable} and {\KoningStabilizable}.
This enables to utilize Lemma \ref{thm:existing_results} {\KoningPiConverge}, {\KoningIHopt}, and {\KoningUnique}.

Using the solution $(\IoptVmat,\IoptFBgain)$,
we set the weight as follows:
\begin{align}
	\forall \StoParam,\;
	{\Weight{}}
	&=
	\dWeight(\StoParam; \riskParam, \estfuncInput(\wildcard;\IoptFBgain)  ,  \estfunccostJ(\wildcard;\IoptVmat)  )
	.\label{eq:problem1_Weight_mod}
\end{align}
The cost function ${\IHWcostJ{\funcInput}{\State_{0}}{\riskParam}}$ with \eqref{eq:problem1_Weight_mod} is equivalent to ${\IHWcostJ{\funcInput}{\State_{0}}{\riskParam}} \SSOCPsetting$
because ${\Expect{\StoParam_{0:\MyT}}}[\prod_{\MybT = 0}^{\MyT}  {\Weight{\MybT}}( \dots )]={\biasE{\StoParam_{0:\MyT}}}[( \dots )]$ holds.
Solving \eqref{eq:def_optFBcontroller} and \eqref{eq:def_mincostJ} reduces to solving the {SSOCP} with ${\biasPDF{\StoParam}}$ in Definition \ref{def:SSOCP}.
From Lemma \ref{thm:existing_results} {\KoningPiConverge} and {\KoningIHopt}, 
the optimal controller and minimum cost are $\estfuncInput(\State;\KoIoptFBgain)$ and $\estfunccostJ(\State;\KoIoptVmat)$, respectively.
Lemma \ref{thm:existing_results} {\KoningUnique} implies the uniqueness of the solution, that is, $\KoIoptVmat=\IoptVmat$ and $\KoIoptFBgain=\IoptFBgain$ hold.
Therefore, 
$\optFBcontroller(\State)=-\IoptFBgain \State=\estfuncInput(\State;\IoptFBgain)$ and $\mincostJ(\State)= \State^{\MyTRANSPO} \IoptVmat \State=\estfunccostJ(\State;\IoptVmat)$
are solutions to \eqref{eq:def_optFBcontroller} and \eqref{eq:def_mincostJ}, respectively, and \eqref{eq:problem1_Weight_mod} is equivalent to \eqref{eq:problem1_Weight}.
This completes the proof.

\section{Proof of Theorem \ref{thm:uniqueness}} \label{pf:uniqueness}

In this proof, a solution to the WSR equations \eqref{eq:def_optVmat_inf} and \eqref{eq:def_optFBgain_inf} for each $\riskParam$ is explicitly denoted by $(\IoptVmat(\riskParam) \succ 0 ,\IoptFBgain(\riskParam) )$, and its vectorization is denoted by ${\optSolPair{\riskParam}}$. 
Note that 
${\ImpTwoF{\SolPair}{0}}=0$ reduces to $\Vmat = {\biasWSRDVmat{\Vmat}}$ and  $\FBgain:={\biasWSRDFBgain{\Vmat}}$ with \eqref{eq:def_optVmat_inf_bias}, \eqref{eq:def_WSRDFBgain_bias}, and ${\biasPDF{\StoParam}}={\PDF{\StoParam}}$.
Because of Assumption \ref{ass:weight} {\AssWeightMSstabilizable} and {\AssWeightZeroSolContained} and Lemma \ref{thm:existing_results} {\KoningPiConverge} and {\KoningUnique}, there exist unique ${\IoptVmat(0)}\succ 0$ and ${\IoptFBgain(0)}={\biasWSRDFBgain{\IoptVmat(0)}}$, that is, ${\optSolPair{0}} \in \DomSolPair$.
We now prove the following statements:
\begin{enumerate}
	
	\item[{\StatementImplicitFunction}] 
	There exist open sets $\DomNeighborriskParam$ and $\DomNeighborSolPair$  
	such that 
	$(0, {\optSolPair{0}} ) \in \DomNeighborriskParam \times \DomNeighborSolPair \subset  \mathbb{R} \times   \mathbb{R}^{(\DimX(\DimX+1)/2) + \DimU \DimX} $ holds and
	there exists a unique $C^{1}$ continuous function ${\optSolPair{\wildcard}} : \DomNeighborriskParam \to \DomNeighborSolPair$ satisfying ${\ImpTwoF{\optSolPair{\riskParam}}{\riskParam}}=0$.
	
	\item[{\StatementBallforImplicitFunction}] 
	There exist $\DomNeighborriskParam$ and $\DomNeighborSolPair$  
	such that
	$\DomNeighborriskParam$ is connected, $\DomNeighborSolPair \subset \DomSolPair$ holds, and  
	all the conditions in {\StatementImplicitFunction} hold.

\end{enumerate}
Firstly, ${\ImpTwoF{\SolPair}{\riskParam}}$ is $C^{2}$ continuous on a neighborhood of $({\optSolPair{0}},0)$ because of Assumption \ref{ass:weight}.
Next, we prove that $\partial_{\SolPair^{\MyTRANSPO}} {\ImpTwoF{\optSolPair{0}}{0}}$ is nonsingular.	
Note that $\VEC{\NotationMat_{1}\NotationMat_{2}\NotationMat_{3}}=(\NotationMat_{3}^{\MyTRANSPO} \otimes \NotationMat_{1}) \VEC{\NotationMat_{2}}$ for given $\NotationMat_{1}$, $\NotationMat_{2}$, and $\NotationMat_{3}$ \citep[Section 3.2.10.2]{Gentle17}. 
Using a relationship similar to $\VECH{ {\Ltrans{}{\FBgain}{\NotationSymMat}{}} } 
= 
{\Elimi{\kronLtrans{\FBgain}}} \VECH{\NotationSymMat}$ in Appendix \ref{sec_supporting},
we derive the several representations of ${\ImpRSR{\SolPair}{\riskParam}}$ and ${\ImpFBgain{\SolPair}{\riskParam}}$:
\begin{align}
	{\ImpRSR{\SolPair}{\riskParam}}
	&
	=
	\Elimi{
		{\WE{\StoParam}{\riskParam}{\FBgain}{\Vmat}}
		[{\DescriptionFBMat{\FBgain}}^{\MyTRANSPO} \otimes	{\DescriptionFBMat{\FBgain}}^{\MyTRANSPO}	
		]	
		-
		{\Identity{\DimX^{2}}}
	}
	\nonumber\\&\quad \times
	\VECH{\Vmat}
	+\VECH{	\FBgain^{\MyTRANSPO} \uCostMat \FBgain + \xCostMat	}
	\nonumber
	\\&
	=
	\VECH[\big]{
		{\WE{\StoParam}{\riskParam}{\FBgain}{\Vmat}} [ \DriftMat^{\MyTRANSPO} \Vmat \DriftMat ]
		-\FBgain^{\MyTRANSPO} {\WE{\StoParam}{\riskParam}{\FBgain}{\Vmat}}[ \InMat^{\MyTRANSPO} \Vmat \DriftMat ]
		\nonumber\\&\quad
		-{\WE{\StoParam}{\riskParam}{\FBgain}{\Vmat}}[ \DriftMat^{\MyTRANSPO} \Vmat \InMat ] \FBgain
		\nonumber\\&\quad
		+\FBgain^{\MyTRANSPO} {\WE{\StoParam}{\riskParam}{\FBgain}{\Vmat}}[ \InMat^{\MyTRANSPO} \Vmat \InMat + \uCostMat ] \FBgain
		+ \xCostMat
		-\Vmat 
	}
	, \label{eq:def_func_imp_2_detail}
	\\
	{\ImpFBgain{\SolPair}{\riskParam}}
	&=
	(
	{\Identity{\DimX}}
	\otimes
	{\WE{\StoParam}{\riskParam}{\FBgain}{\Vmat}} [ \InMat^{\MyTRANSPO} \Vmat \InMat  +  \uCostMat]
	)
	\VEC{\FBgain}
	\nonumber\\&\quad
	-
	\VEC{	{\WE{\StoParam}{\riskParam}{\FBgain}{\Vmat}} [ \InMat^{\MyTRANSPO} \Vmat \DriftMat] }
	.\label{eq:def_func_imp_1_detail}
\end{align}
Assumption \ref{ass:desired_weight} with $\riskParam=0$ implies that for every $(\StoParam, \estfuncInput(\wildcard;\FBgain)  ,  \estfunccostJ(\wildcard;\Vmat))$, we have $\dWeight(\StoParam; \riskParam, \estfuncInput(\wildcard;\FBgain)  ,  \estfunccostJ(\wildcard;\Vmat)  )=1$  and thus $\partial \dWeight(\StoParam; \riskParam, \estfuncInput(\wildcard;\FBgain)  ,  \estfunccostJ(\wildcard;\Vmat)  )/\partial \SolPair=0$.
The following partial derivatives are obtained from \eqref{eq:def_func_imp_2_detail} and \eqref{eq:def_func_imp_1_detail}:
\begin{align}
	&\frac{\partial }{\partial \VECH{\Vmat}^{\MyTRANSPO}}
	{\ImpRSR{\optSolPair{0}}{0}}
	\nonumber\\&
	=
	\Elimi{
		{\Expect{\StoParam}}
		[{\DescriptionFBMat{{\IoptFBgain(0)}}}^{\MyTRANSPO} 
		\otimes	{\DescriptionFBMat{{\IoptFBgain(0)}}}^{\MyTRANSPO}	
		]	
	}
		\nonumber\\&\quad -	{\Identity{\DimX(\DimX+1)/2}}
	,\label{eq:def_func_imp_2_partial_Pi_theta0}
	\\
	&\frac{\partial }{\partial \VEC{\FBgain}^{\MyTRANSPO} }
	{\ImpFBgain{\optSolPair{0}}{0}}
	=
	{\Identity{\DimX}}
	\otimes
	{\Expect{\StoParam}}[ \InMat^{\MyTRANSPO} {\IoptVmat(0)} \InMat  +  \uCostMat]
	.\label{eq:def_func_imp_1_partial_L_theta0}
\end{align} 
Note that
$\NotationSymMat 
\mapsto 
{\Expect{\StoParam}}
[{\DescriptionFBMat{{\IoptFBgain(0)}}}^{\MyTRANSPO} 
\NotationSymMat
{\DescriptionFBMat{{\IoptFBgain(0)}}}
]	$
is stable by Lemma \ref{thm:existing_results} {\KoningUnique}.
Because of Definition \ref{def:stable_trans},
the spectral radius of
$	\Elimi{
			{\Expect{\StoParam}}
			[{\DescriptionFBMat{{\IoptFBgain(0)}}}^{\MyTRANSPO} 
			\otimes
			{\DescriptionFBMat{{\IoptFBgain(0)}}}^{\MyTRANSPO}	
			]	
		}
$, that is, maximum absolute value of the eigenvalues, is less than 1 and thus \eqref{eq:def_func_imp_2_partial_Pi_theta0} is nonsingular.
Because ${\Expect{\StoParam}}[ \InMat^{\MyTRANSPO} {\IoptVmat(0)} \InMat  +  \uCostMat]\succ 0$, the block diagonal matrix $\partial_{\VEC{\FBgain}^{\MyTRANSPO}} {\ImpFBgain{\optSolPair{0}}{0}}$ in (\ref{eq:def_func_imp_1_partial_L_theta0}) is nonsingular. 
The following partial derivative is calculated from (\ref{eq:def_func_imp_2_detail}):
\begin{align}
	&\frac{\partial}{\partial {\El{\FBgain}{\IDEl,\IDbEl}}}
	{\ImpRSR{\SolPair}{0}}
	\nonumber\\&
	=
	\VECH[\Big]{
		\frac{\partial \FBgain^{\MyTRANSPO}}{\partial {\El{\FBgain}{\IDEl,\IDbEl}}} 
		\Big(
		{\Expect{\StoParam}}[ \InMat^{\MyTRANSPO} \Vmat \InMat  + \uCostMat ]\FBgain
		- {\Expect{\StoParam}}[ \InMat^{\MyTRANSPO} \Vmat \DriftMat ]
		\Big)
		\nonumber\\&\quad		
		+
		\Big(
		\FBgain^{\MyTRANSPO} {\Expect{\StoParam}}[ \InMat^{\MyTRANSPO} \Vmat \InMat  + \uCostMat ]
		- {\Expect{\StoParam}}[ \DriftMat^{\MyTRANSPO} \Vmat \InMat  ]
		\Big)
		\frac{\partial \FBgain}{\partial {\El{\FBgain}{\IDEl,\IDbEl}}} 
	} 
	.\label{eq:def_func_imp_2_partial_L_theta0}		
\end{align} 
We have
$\partial_{\VEC{\FBgain}^{\MyTRANSPO}} {\ImpRSR{\optSolPair{0}}{0}}=0$ 
by substituting ${\optSolPair{0}}$ with ${\IoptFBgain(0)}=
{\Expect{\StoParam}}[ \InMat^{\MyTRANSPO} \IoptVmat(0) \InMat  + \uCostMat ]^{-1}{\Expect{\StoParam}}[ \InMat^{\MyTRANSPO} \IoptVmat(0) \DriftMat ]$ into \eqref{eq:def_func_imp_2_partial_L_theta0}. 
Thus, $\partial_{\SolPair^{\MyTRANSPO}} {\ImpTwoF{\optSolPair{0}}{0}}$ is the block triangular matrix that  is nonsingular because we have nonzero determinants %
$\det|\partial_{\VEC{\Vmat}^{\MyTRANSPO}} {\ImpRSR{\optSolPair{0}}{0}}|\neq0$ and %
$\det|\partial_{\VEC{\FBgain}^{\MyTRANSPO}} {\ImpFBgain{\optSolPair{0}}{0}}|\neq0$ \citep[Section 3.1.9.7]{Gentle17}. %
Then, {\StatementImplicitFunction} holds according to the implicit function theorem \citep[Theorem 5]{Oliveira2013}.
We obtain {\StatementBallforImplicitFunction} by choosing $\DomNeighborriskParam$ as a small ball of $0$ because $\DomNeighborriskParam$ and $\DomNeighborSolPair$ contain open balls of $0$ and ${\optSolPair{0}}$, respectively, and ${\optSolPair{\wildcard}}$ is continuous.

Using {\StatementBallforImplicitFunction}, we prove Theorem \ref{thm:uniqueness} as follows. 
Let $\DomBSolPair:=\DomSolPair \setminus \DomNeighborSolPair $.
Because $\DomBSolPair$ is a bounded closed set,
$\|{\ImpTwoF{\SolPair}{0}}\|$ is continuous on $\DomBSolPair$,
and the uniqueness of ${\optSolPair{0}}$ indicates ${\ImpTwoF{\SolPair}{0}} \neq 0$ on $\DomBSolPair$,
there exists $\BoundImpTwoF>0$ that satisfies $\|{\ImpTwoF{\SolPair}{0}}\|>\BoundImpTwoF$ for every $\SolPair \in \DomBSolPair$. 
For some $\tmpLBriskParam<0$ and $\tmpUBriskParam>0$ that satisfy $[\tmpLBriskParam,\tmpUBriskParam] \subset \DomNeighborriskParam$,
$\partial_{\riskParam} {\ImpTwoF{\SolPair}{\riskParam}}$ is bounded on the bounded closed set $\DomBSolPair \times [\tmpLBriskParam,\tmpUBriskParam]$ because of the continuity.
Namely, there exists $\BoundpartialImpTwoF>0$ such that $\| \partial_{\riskParam} {\ImpTwoF{\SolPair}{\riskParam}}\| \leq \BoundpartialImpTwoF$.
Here, we set $\uniqueLBriskParam = \max \{   -\BoundImpTwoF/ \BoundpartialImpTwoF , \tmpLBriskParam  \}<0$ and 
$\uniqueUBriskParam = \min \{   \BoundImpTwoF/ \BoundpartialImpTwoF , \tmpUBriskParam  \}>0$.
Then, for any $(\SolPair,\riskParam) \in \DomBSolPair \times [\uniqueLBriskParam, \uniqueUBriskParam]$, we have 
\begin{align}	
	\|{\ImpTwoF{\SolPair}{\riskParam}}\|
	\geq
	\|{\ImpTwoF{\SolPair}{0}}\| -  \BoundpartialImpTwoF |\riskParam| 
	>
	 \BoundImpTwoF -  \BoundpartialImpTwoF |\riskParam|     \geq 0
	.
\end{align}	
Thus, there exists no solution to ${\ImpTwoF{\SolPair}{\riskParam}}=0$ on  $(\DomSolPair \setminus \DomNeighborSolPair ) \times [\uniqueLBriskParam, \uniqueUBriskParam]$.
Meanwhile, there exist a unique solution $({\optSolPair{\riskParam}},\riskParam)$ on  $\DomNeighborSolPair \times [\uniqueLBriskParam, \uniqueUBriskParam]$ because  $[\uniqueLBriskParam, \uniqueUBriskParam] \subset \DomNeighborriskParam$ holds.
Therefore, 
there exists a unique ${\optSolPair{\wildcard}} : [\uniqueLBriskParam, \uniqueUBriskParam] \to \DomSolPair$ that satisfies ${\ImpTwoF{\optSolPair{\riskParam}}{\riskParam}}=0$ and $C^{1}$ continuous on $(\uniqueLBriskParam, \uniqueUBriskParam)$.
The positive definiteness $\IoptVmat \succ 0$ holds form the definition of $\DomVmat$.
This completes the proof.

\section{Proof of Theorem \ref{thm:solution_WSR}} \label{pf:solution_WSR}

We prove the statement by taking the limit $\MybT \rightarrow \infty$ in a manner similar to \citep[Theorem 5.1]{Koning82}.	
Because
${\WSRDVFB{\wildcard}{\wildcard}{\riskParam}}:=[ {\WSRDVmat{\wildcard}{\wildcard}{\riskParam}},{\WSRDFBgain{\wildcard}{\wildcard}{\riskParam}}{}^{\MyTRANSPO}]$ is continuous at $(\IestVmat,\IestFBgain)$,
using  \eqref{eq:lim_WSRDeq} yields
$\lim_{\MybT \rightarrow \infty} 
{\WSRDVFB{\optVmat{\MybT}}{\optFBgain{\MybT}}{\riskParam}}
=
{\WSRDVFB{\IestVmat}{\IestFBgain}{\riskParam}}
$ as follows:
\begin{align}	
	\forall & \MyEpsilonForConv>0,  \quad	
	\exists \MyDeltaForConv>0, \quad
	\exists \MybT^{\prime}>0,
	\nonumber \\
	\MybT > \MybT^{\prime}
	&
	\Rightarrow
	\FrobeniusNorm{ 
		[{\optVmat{\MybT}},{\optFBgain{\MybT}}{}^{\MyTRANSPO}]
		- 
		[\IestVmat	,		\IestFBgain{}^{\MyTRANSPO}]	 
	} < \MyDeltaForConv
	\nonumber \\&
	\Rightarrow 
	\FrobeniusNorm{ 
		{\WSRDVFB{\optVmat{\MybT}}{\optFBgain{\MybT}}{\riskParam}}
		-
		{\WSRDVFB{\IestVmat}{\IestFBgain}{\riskParam}}
	}<\MyEpsilonForConv
	,
\end{align}	
where $\FrobeniusNorm{\wildcard}$ is the Frobenius norm.
Thus, we obtain
$[\IestVmat,\IestFBgain{}^{\MyTRANSPO}]
=
\lim_{\MybT \rightarrow \infty} 
[{\optVmat{\MybT}},{\optFBgain{\MybT}}{}^{\MyTRANSPO}]
=
\lim_{\MybT \rightarrow \infty} 
[{\optVmat{\MybT+1}},{\optFBgain{\MybT+1}}{}^{\MyTRANSPO}]
$ $=
{\WSRDVFB{\IestVmat}{\IestFBgain}{\riskParam}}
$.
Next, \eqref{eq:def_WSRD} is transformed into 
\begin{align}
	{\optVmat{\MybT+1}}
	&=
	{\WE{\StoParam}{\riskParam}{\optFBgain{\MybT}}{\optVmat{\MybT}}}
	[(\DriftMat -\InMat {\WSRDFBgain{\optVmat{\MybT}}{\optFBgain{\MybT}}{\riskParam}} )^{\MyTRANSPO} 
	\nonumber\\&\quad \times
	{\optVmat{\MybT}} 
	(\DriftMat -\InMat {\WSRDFBgain{\optVmat{\MybT}}{\optFBgain{\MybT}}{\riskParam}} )
	]
	\nonumber\\&\quad
	+  
	{\WSRDFBgain{\optVmat{\MybT}}{\optFBgain{\MybT}}{\riskParam}}^{\MyTRANSPO} 
	\uCostMat 
	{\WSRDFBgain{\optVmat{\MybT}}{\optFBgain{\MybT}}{\riskParam}}   
	+ \xCostMat
	.
\end{align}
Because $\dWeight(\StoParam; \riskParam, \estfuncInput(\wildcard;\FBgain)  ,  \estfunccostJ(\wildcard;\Vmat)  )\geq 0$ holds,
the positive semidefiniteness of $\IiniVmat \succeq 0$ implies ${\optVmat{\MybT}}\succeq \xCostMat \succ 0$ for any $\MybT$, which yields $\IestVmat\succ 0$.
This completes the proof.

\section{Proof of Theorem \ref{thm:MyNewton}} \label{pf:MyNewton}

In this proof, let $\| \NotationMat \|$ be the induced norm of a given matrix $\NotationMat$.
Because of Assumption \ref{ass:weight} {\AssWeightSmooth}, ${\ImpTwoF{\SolPair}{\riskParam}}$ is $C^{2}$ continuous, and thus $\partial_{\SolPair^{\MyTRANSPO}}{\ImpTwoF{\SolPair}{\riskParam}}$ is Lipschitz continuous on $\DomSolPair \times[\smoothLBriskParam,\smoothUBriskParam] $.
Let $\NewtonLipschitz$ be a Lipschitz constant of $\partial_{\SolPair^{\MyTRANSPO}}{\ImpTwoF{\SolPair}{\riskParam}}$ on this set.
In addition, for some $ \uniqueLBriskParam\geq \smoothLBriskParam$ and $\uniqueUBriskParam\leq \smoothUBriskParam$, for every $\riskParam \in [\uniqueLBriskParam,\uniqueUBriskParam]$, ${\optSolPair{\riskParam}}$ is a unique solution by Theorem \ref{thm:uniqueness}.
Based on the proofs in {\citep[Lemma 4.3.1 and Theorem 5.1.1]{Kelley95}},
the q-quadratical convergence holds in Lemma \ref{def:original_Newton} if  the standard assumptions and ${\iteSolPair{0}} = {\optSolPair{0}} \in {\BallSolPair{\reqBallRadius(\riskParam)}{\riskParam}}$ hold for a scalar ${\reqBallRadius(\riskParam)}$ that satisfies
$ \NewtonLipschitz \| \partial_{\SolPair^{\MyTRANSPO}}{\ImpTwoF{\optSolPair{\riskParam}}{\riskParam}}^{-1} \|  \reqBallRadius(\riskParam) < 1/2$.
We define 
$\reqBallRadius(\riskParam):=\NewtonScalar/(2 \NewtonLipschitz \| \partial_{\SolPair^{\MyTRANSPO}}{\ImpTwoF{\optSolPair{\riskParam}}{\riskParam}}^{-1} \| )$ with a constant $\NewtonScalar \in (0,1)$.
Thus,
Theorem \ref{thm:MyNewton} is derived if the following statements hold.
\begin{enumerate}
	\item[{\StatementSetStandardAss}] 
	There exist an open set $\DomSaSolPair \subset \DomSolPair$ containing ${\optSolPair{0}}$,
	$\NewtonUBriskParam>0$, and $\NewtonLBriskParam < 0$ such that for every $\riskParam \in (\NewtonLBriskParam,\NewtonUBriskParam)$,
	the standard assumptions on $(\DomSaSolPair,\riskParam)$ hold.

	\item[{\StatementInitialGuess}]
	For any open set $\DomSaSolPair \subset \DomSolPair$ containing ${\optSolPair{0}}$,
	there exist $\DefBallRadius>0$, $\NewtonLBriskParam<0$, and $\NewtonUBriskParam>0$ such that
	for every $\riskParam \in (\NewtonLBriskParam,\NewtonUBriskParam)$,
	we have ${\BallSolPair{\DefBallRadius}{\riskParam}} \subseteq {\BallSolPair{\reqBallRadius(\riskParam) }{\riskParam}}$ and ${\iteSolPair{0}}={\optSolPair{0}} \in {\BallSolPair{\DefBallRadius}{\riskParam}} \subset \DomSaSolPair$.
	
\end{enumerate}

We prove that the standard assumption {\StandardAssA} holds for deriving {\StatementSetStandardAss}.
According to Theorem \ref{thm:uniqueness},  
there exists a unique solution ${\optSolPair{\riskParam}} \in \DomSolPair$ that is    $C^{1}$ continuous in $\riskParam$ on $(\uniqueLBriskParam,\uniqueUBriskParam)$.  
Let $\DomSaSolPair\subset \DomSolPair$ be any open ball with the center ${\optSolPair{0}}$.
There exist  $\NewtonLBriskParam \in [\uniqueLBriskParam,0)$ and $\NewtonUBriskParam \in (0,\uniqueUBriskParam]$ such that ${\optSolPair{\riskParam}} \in \DomSaSolPair$ holds for any $\riskParam \in (\NewtonLBriskParam,\NewtonUBriskParam)$ because of the continuity.
This indicates the standard assumption {\StandardAssA} for any open ball $\DomSaSolPair$ and $\riskParam \in (\NewtonLBriskParam,\NewtonUBriskParam) \subseteq (\uniqueLBriskParam,\uniqueUBriskParam)$.

Next, we prove that the standard assumption {\StandardAssB} holds.
Because $\partial_{\SolPair^{\MyTRANSPO}}{\ImpTwoF{\SolPair}{\riskParam}}$ is Lipschitz continuous on $\DomSolPair \times[\smoothLBriskParam,\smoothUBriskParam] $ according to Assumption \ref{ass:weight},
the standard assumption {\StandardAssB} holds for any open set $\DomSaSolPair \subset \DomSolPair$ containing ${\optSolPair{0}}$ and any $\riskParam \in [\smoothLBriskParam,\smoothUBriskParam]$.

We prove that the standard assumption {\StandardAssC} holds.
According to the standard assumption {\StandardAssB}, 
let ${\Distance{\SolPair}{\riskParam}} := [\SolPair^{\MyTRANSPO}, \riskParam]-[{\optSolPair{0}}^{\MyTRANSPO},0]$, 
and we have
$	\|\partial_{\SolPair^{\MyTRANSPO}}{\ImpTwoF{\SolPair}{\riskParam}}
-
\partial_{\SolPair^{\MyTRANSPO}}{\ImpTwoF{\optSolPair{0}}{0}}\|
\leq
\NewtonLipschitz 	\|{\Distance{\SolPair}{\riskParam}}\|$,
where $\partial_{\SolPair^{\MyTRANSPO}}{\ImpTwoF{\optSolPair{0}}{0}}$ is nonsingular because of  the proof of Theorem \ref{thm:uniqueness}.
We obtain the following result based on \citep[the proof of Lemma 4.3.1]{Kelley95}: 	
\begin{align}	
	&
	\|
	{\Identity{\DimX(\DimU+(\DimX+1)/2)}}-
	\partial_{\SolPair^{\MyTRANSPO}}{\ImpTwoF{\optSolPair{0}}{0}}^{-1}
	\partial_{\SolPair^{\MyTRANSPO}}{\ImpTwoF{\SolPair}{\riskParam}}
	\|
	\nonumber\\&
	=
	\|
	\partial_{\SolPair^{\MyTRANSPO}}{\ImpTwoF{\optSolPair{0}}{0}}^{-1}
	(
	\partial_{\SolPair^{\MyTRANSPO}}{\ImpTwoF{\optSolPair{0}}{0}}
	-
	\partial_{\SolPair^{\MyTRANSPO}}{\ImpTwoF{\SolPair}{\riskParam}}
	)
	\|
	\nonumber\\&
	\leq
	\NewtonLipschitz
	\|\partial_{\SolPair^{\MyTRANSPO}}{\ImpTwoF{\optSolPair{0}}{0}}^{-1}\|
	\|{\Distance{\SolPair}{\riskParam}}\|
	.
\end{align}	
If the right hand side of this inequality is less than 1,
the Banach Lemma \citep[Theorem 1.2.1]{Kelley95}	implies that $\partial_{\SolPair^{\MyTRANSPO}}{\ImpTwoF{\SolPair}{\riskParam}}$ is nonsingular.
In other words, this nonsingular property holds on some open ball $\BallNonsingular \subset\DomSolPair \times[\smoothLBriskParam,\smoothUBriskParam] $ with the center $({\optSolPair{0}},{0})$ and a radius less than
$
1/(
\NewtonLipschitz
\|\partial_{\SolPair^{\MyTRANSPO}}{\ImpTwoF{\optSolPair{0}}{0}}^{-1}\|
)
$.
There exist an open set $\DomSaSolPair  \subset \DomSolPair$, $\NewtonLBriskParam<0$, and $\NewtonUBriskParam>0$ such that
$\DomSaSolPair\times  (\NewtonLBriskParam,\NewtonUBriskParam) \subset \BallNonsingular$  and
${\optSolPair{\riskParam}} \in \DomSaSolPair$  for every $\riskParam \in (\NewtonLBriskParam,\NewtonUBriskParam)$ are satisfied because of the continuity of ${\optSolPair{\riskParam}}$.
Subsequently, for every $\riskParam \in (\NewtonLBriskParam,\NewtonUBriskParam)$,
we have $( {\optSolPair{\riskParam}}, \riskParam) \in \BallNonsingular$; thus, the standard assumptions  {\StandardAssC} on $(\DomSaSolPair,\riskParam)$ hold.
From these results, {\StatementSetStandardAss} holds.

Next, 
we prove the statement {\StatementInitialGuess}.
Recall that $\partial_{\SolPair^{\MyTRANSPO}}{\ImpTwoF{\SolPair}{\riskParam}}$ is nonsingular on the ball  $\BallNonsingular$ with the center $({\optSolPair{0}},{0})$.
Because 
${\optSolPair{\riskParam}}$ is continuous,
$\|\partial_{\SolPair^{\MyTRANSPO}}{\ImpTwoF{\optSolPair{\riskParam}}{\riskParam}}^{-1}\|$ is continuous and thus bounded on a neighborhood of $\riskParam=0$.
Using the boundedness, for some $\NewtonLLBriskParam<0$ and $\NewtonUUBriskParam>0$, there exists $\LBreqBallRadius:=\inf_{ \riskParam \in (\NewtonLLBriskParam,\NewtonUUBriskParam) } \reqBallRadius(\riskParam) >0$.
In addition,
for any open set $\DomSaSolPair \subset \DomSolPair$ containing ${\optSolPair{0}}$,
there exists $\DefBallRadius\in (0,\LBreqBallRadius)$ such that  $ {\optSolPair{0}} \in {\BallSolPair{\DefBallRadius}{0}} \subset \DomSaSolPair$ holds.
Therefore,
by using the continuity of $\optSolPair{\riskParam}$,
for such $\DefBallRadius\in (0,\LBreqBallRadius)$, there exist $\NewtonLBriskParam\in [\NewtonLLBriskParam,0)$ and $\NewtonUBriskParam \in (0,\NewtonUUBriskParam]$ such that
for every $\riskParam \in (\NewtonLBriskParam,\NewtonUBriskParam)$,
we have $ {\optSolPair{0}} \in {\BallSolPair{\DefBallRadius}{\riskParam}} \subseteq {\BallSolPair{\LBreqBallRadius}{\riskParam}} \subseteq {\BallSolPair{\reqBallRadius(\riskParam) }{\riskParam}}$ and $ {\BallSolPair{\DefBallRadius}{\riskParam}} \subset \DomSaSolPair$.
This indicates that {\StatementInitialGuess} holds.
This completes the proof.

\section{Proof of Theorem \ref{thm:ms_stable_existence}} \label{pf:ms_stable_existence}

We can prove this theorem by using Theorem \ref{thm:uniqueness} and Lemma \ref{thm:existing_results} with ${\biasPDF{\StoParam}}={\PDF{\StoParam}}$.
Because the feedback system \eqref{eq:def_FBsys} is MS stabilizable by Assumption \ref{ass:weight}, there exist unique ${\IoptVmat(0)}\succ 0$ and ${\IoptFBgain(0)}={\biasWSRDFBgain{\IoptVmat(0)}}$ because of Lemma \ref{thm:existing_results} {\KoningPiConverge} and {\KoningUnique}.
Because Theorem \ref{thm:uniqueness} states that ${\IoptFBgain(\riskParam)}$ is continuous in $\riskParam \in (\uniqueLBriskParam,\uniqueUBriskParam)$, 
the matrix $\QuadFBMat(\riskParam):=
{\Elimi{
{\Expect{\StoParam}}[{\DescriptionFBMat{\IoptFBgain(\riskParam)}} \otimes {\DescriptionFBMat{\IoptFBgain(\riskParam)}} ]
}} \in \mathbb{R}^{\DimexX \times \DimexX}$ is also continuous on $(\uniqueLBriskParam,\uniqueUBriskParam)$, where $\DimexX:=\DimX(\DimX+1)/2$.
There exist continuous functions $({\EigQuadFBMat{1}{\riskParam}},{\EigQuadFBMat{2}{\riskParam}}, \dots, {\EigQuadFBMat{\DimexX}{\riskParam}}  )$ on $(\uniqueLBriskParam,\uniqueUBriskParam)$ such that their values are equal to the repeated eigenvalues  of $\QuadFBMat(\riskParam)$ \citep[Theorem 5.2]{Kato84}.
The spectral radius $\max_{\IDEl} | {\EigQuadFBMat{\IDEl}{\riskParam}} |$ is also continuous.
Lemma \ref{thm:existing_results} {\KoningUnique} indicates that  
${\NonArgLtrans{\IoptFBgain(0)}{}}$  with ${\biasPDF{\StoParam}}={\PDF{\StoParam}}$ is stable.
The spectral radius $\max_{\IDEl} | {\EigQuadFBMat{\IDEl}{0}} |$ is less than 1 according to Definition \ref{def:stable_trans}.
Because of the continuity of $\max_{\IDEl} | {\EigQuadFBMat{\IDEl}{\riskParam}} |$ on $(\uniqueLBriskParam,\uniqueUBriskParam)\ni 0$,
there exists  $\StableUBriskParam > 0$ and $\StableLBriskParam < 0$ such that for every $\riskParam \in ( \StableLBriskParam , \StableUBriskParam)$, we have $\max_{\IDEl} | {\EigQuadFBMat{\IDEl}{\riskParam}} |<1$ and thus ${\NonArgLtrans{\IoptFBgain(\riskParam)}{}}$  with ${\biasPDF{\StoParam}}={\PDF{\StoParam}}$ is stable  by Definition \ref{def:stable_trans}.
Then,  the system \eqref{eq:def_FBsys} with $\funcInput(\State)=-{\IoptFBgain(\riskParam)} \State$ is MS stable by Lemma \ref{thm:existing_results} {\KoningStable}.
This completes the proof.

\section{Proof of Theorem \ref{thm:wms_stable}} \label{pf:wms_stable}

In a manner similar to the proof of Theorem \ref{thm:solution_problem_1},
let us consider the system \eqref{eq:def_sys} with $\StoParam_{\MyT}$ that obeys the biased PDF ${\biasPDF{\StoParam}}
=
\dWeight(\StoParam_{\MybT}; \riskParam, \optFBcontroller ,  \mincostJ ) 
\PDF{\StoParam}$ called the biased system.
Subsequently, the WSR equations  \eqref{eq:def_optVmat_inf} and \eqref{eq:def_optFBgain_inf} reduce to
$\Vmat = {\biasWSRDVmat{\Vmat}}$ and  $\FBgain:={\biasWSRDFBgain{\Vmat}}$ with 
\eqref{eq:def_optVmat_inf_bias} and \eqref{eq:def_WSRDFBgain_bias}, respectively.
Because of $\IoptVmat\succ 0$ is a solution to $\Vmat = {\biasWSRDVmat{\Vmat}}$,
we have 
${\biasE{\StoParam_{0:\MyT}}}[\State_{\MyT}^{\MyTRANSPO} {\NotationSymMat} \State_{\MyT}]	
=\State_{0}^{\MyTRANSPO} {\Ltrans{}{\IoptFBgain}{\NotationSymMat}{\MyT}} \State_{0}$	
with a stable $ {\Ltrans{}{\IoptFBgain}{\wildcard}{}} $
from Lemma \ref{thm:existing_results} {\KoningTransStable} and {\KoningTransition}.
Lemma \ref{thm:def_stable2} implies that $
{\Expect{\StoParam_{0:\MyT}}}[ {\WmsWeight{\MyT}} \|\State_{\MyT}\|^{2} ]
=
{\Expect{\StoParam_{0:\MyT}}}[ \prod_{\MybT=0}^{\MyT}
\dWeight(\StoParam_{\MybT}; \riskParam, \optFBcontroller ,  \mincostJ )  \|\State_{\MyT}\|^{2} ]
=
{\biasE{\StoParam_{0:\MyT}}}[  \|\State_{\MyT}\|^{2} ]
=
\State_{0}^{\MyTRANSPO} {\Ltrans{}{\IoptFBgain}{\Identity{\DimX}}{\MyT}}  \State_{0}
\to 0$ as $\MyT \to \infty$ in a manner similar to Lemma \ref{thm:existing_results} {\KoningStable}.
This completes the proof.

\end{document}